\documentclass[10pt,a4paper]{article}
\usepackage[utf8]{inputenc}
\usepackage[english]{babel}
\usepackage[T1]{fontenc}
\usepackage{amsmath}
\usepackage{amsfonts}
\usepackage{graphicx}
\usepackage{array}
\usepackage[left=5cm,right=5cm,top=4cm,bottom=4cm]{geometry}

\usepackage{adjustbox}
\usepackage{transparent}

\usepackage{isomath}
\usepackage{mathtools}
\usepackage{amsthm}
\usepackage{slashed}
\usepackage{amssymb}
\usepackage{enumitem}
\usepackage{hyperref}
\hypersetup{colorlinks=true, 
citecolor=red}
\usepackage[new]{old-arrows}

\theoremstyle{definition}
\newtheorem{mdef}{{Definition}}[section]

\theoremstyle{definition}
\newtheorem{mex}{Example}[section]

\theoremstyle{definition}
\newtheorem{mrmk}{{Remark}}[section]

\theoremstyle{plain}
\newtheorem{mth}{Theorem}[section]

\theoremstyle{plain}
\newtheorem{mlem}{{Lemma}}[section]

\theoremstyle{plain}
\newtheorem{mprop}{{Proposition}}[section]

\theoremstyle{plain}
\newtheorem{mcor}{{Corollary}}[section]

\usepackage[symbol]{footmisc}

\newcommand*\samethanks[1][\value{footnote}]{\footnotemark[#1]}
    
\makeatletter
\def\keywords{\xdef\@thefnmark{}\@footnotetext}
\makeatother 

\DeclareMathOperator{\Ima}{Im}

\DeclarePairedDelimiter{\abs}{\lvert}{\rvert}
\makeatletter
\let\oldabs\abs
\def\abs{\@ifstar{\oldabs}{\oldabs*}}
\makeatother

\newcommand{\rest}[1]{\left.{#1}\right\vert}

\newcommand{\bigslant}[2]{{\raisebox{.2em}{$#1$}\left/\raisebox{-.2em}{$#2$}\right.}}

\def\XXint#1#2#3{{\setbox0=\hbox{$#1{#2#3}{\int}$}
     \vcenter{\hbox{$#2#3$}}\kern-.5\wd0}}

\usepackage{tikz-cd} 
\usetikzlibrary{positioning, shapes, arrows.meta}
\tikzset{
    answer/.style={rectangle, draw, text width=15em, text badly centered, node distance=1cm, inner sep=0pt, minimum height=4em},
    block/.style={rectangle, draw, text width=10em, text centered},
     block2/.style={rectangle, draw, text width=5em, text centered},
     block3/.style={rectangle, draw, text width=5em, text centered, color=white},
}

\DeclareMathSymbol{\upLambda}{\mathalpha}{operators}{3}

\usepackage{pgfplots}
\usetikzlibrary{backgrounds}
\pgfplotsset{compat=1.14}
\usepackage{mathrsfs}
\usetikzlibrary{arrows}

\title{Noncommutative Differential Geometry on Infinitesimal Spaces}
\author{ Damien Tageddine\footnote{Department of Mathematics and Statistics, McGill University, Montreal, Canada} \thanks{Email address: damien.tageddine@mail.mcgill.ca}
\and
 Jean-Christophe Nave \samethanks[1] \thanks{Corresponding author. Email address: jean-christophe.nave@mcgill.ca} }
\date{}
\begin{document}
\maketitle
\begin{abstract}
In this paper, we use the language of noncommutative differential geometry to formalise discrete differential calculus. We begin with a brief review of inverse limit of posets as an approximation of topological spaces. We then show how to associate a $C^*$-algebra over a poset, giving it a piecewise-linear structure. Furthermore, we explain how dually the algebra of continuous function $C(M)$ over a manifold $M$ can be approximated by a direct limit of $C^*$-algebras over posets. Finally, in the spirit of noncommutative differential geometry, we define a finite dimensional spectral triple on each poset. We show how the usual finite difference calculus is recovered as the eigenvalues of the commutator with the Dirac operator. We prove a convergence result in the case of the $d$-lattice in $\mathbb{R}^d$ and for the torus $\mathbb{T}^d$.
\end{abstract}

\keywords{2020 \emph{Mathematics Subject Classification.} 65N99 58B34}
\keywords{\emph{Key words and phrases.} Noncommutative differential geometry, discrete differential calculus, finite difference, C*-algebras, Dirac operators.}
    
\newpage

{
  \hypersetup{linkcolor=black}
  \tableofcontents
}
\newpage
\section{Introduction}
The general motivation for the present work is the discretization of partial differential equations (PDE). This paper aims at laying down the foundation of a broad framework to study discrete differential calculus in a discretization-free fashion. Using the tools of noncommutative differential geometry, we establish a geometric formalism of finite difference calculus in order to tackle the problem of differential operators approximations. We start by recalling general results on approximation of a compact Hausdorff space $M$ by a sequence of ordered simplicial complexes (Proposition \ref{main0}). We then show that the space of continuous functions $C(M)$ can be replaced by a sequence of $C^*$-algebras over each simplicial complex (Theorem \ref{main1} and Proposition \ref{main2}). Finally, after introducing a differential structure on these $C^*$-algebras, we show that the usual finite difference approximations are recovered as eigenvalues of the exterior derivative operator (Proposition \ref{main5}). The convergence of this differential operator to the classical de Rham differential is shown in the case of the $n$-dimensional lattice (Proposition \ref{spectral}) and similarly in the case of the $n$-dimensional torus.
\paragraph*{Related approaches and background} The approximation theory of partial differential equations (PDE) can take several aspects. The various methods rely on the intuitive geometric idea that the fine structure of a space $M$ (one can think of a domain in $\mathbb{R}^n$ or a smooth manifold) is discrete. The resulting discretized space, say $X$, is governed by a parameter — being a grid spacing, the size of a mesh or a time step for example — denoted by $h$, $\varepsilon$ or $\Delta x$, which plays the role of an infinitesimal. In the rest of this work, we will loosely call this type of discrete space \textit{infinitesimal space}. Information extracted from the continuous space can be represented by a family of morphisms $(\chi_x)_{x\in X}$ with
\begin{equation*}
\chi_x:C^\infty(M)\rightarrow \mathbb{C}, \quad \chi_x(f)=f(x),
\end{equation*}
which can be related to either sampling morphisms in finite difference (volume) language, or nodal basis in finite element denominations. These maps encapsulate the local data available from the algebra of functions over the continuous space $M$.\par
The geometric approach of discrete differential calculus has been pioneered by Whitney in his work on geometric integration theory \cite{whitney_geometric_1957}. The classical differential forms can be interpreted as cochains when restricted to a simplicial complex $K$ by means of the de Rham map:
\begin{equation*}
C:\Omega^p(M)\rightarrow C^p(K,\mathbb{Z}), \quad C(\omega):= \sigma \mapsto \left\langle \omega, \sigma \right\rangle.
\end{equation*}
Vice-versa, a cochain can be used to define a differential form using Whitney's interpolation map $\mathcal{W}:C^p(K,\mathbb{Z})\rightarrow \Omega^p(M)$,
\begin{equation*}
\mathcal{W}(x_0,\dots,x_p)=p!\sum_{i=0}^p(-1)^i\lambda_id\lambda_0\wedge \cdots \wedge \widehat{d\lambda_i}\wedge \cdots\wedge d\lambda_p.
\end{equation*}
This viewpoint has then been successfully used in lattice (quantum) field theory in \cite{wilson_confinement_1974, sen_geometric_2000, adams_r-torsion_1996} and in computational electromagnetism \cite{bossavit_generalized_2001,stern_geometric_2015}.\\
Moreover, the idea of deriving a discrete theory that parallels the continuous one has then been further explored by Hirani in the discrete exterior calculus (DEC) \cite{hirani_discrete_2003} and subsequently developed by Desbrun et al. \cite{desbrun_discrete_2005}. In DEC, the point of view — which is also shared to some extent by our work — is that the discrete theory can, and indeed should, stand on its own right. The authors base their approach on simplicial complexes and its differential calculus on chains and cochains. In that setting, a differential form is an element in the dual of the space of chains. The basic data in the theory is given  by the triple $(K,\Omega^*(K),d)$ where $K$ is a simplicial complex, $d$ is the coboundary map and $\Omega^*(K)$ the space of cochains. To this, one adds a Hodge-star map:
\begin{equation*}
        (K,\Omega^*(K),d), \quad *:\Omega^k(K)\rightarrow \Omega^{d-k}(*K)
\end{equation*}
where $*K$ is the dual simplicial complex.\par
In the realm of finite element method, the pioneering work of Arnold et al. \cite{arnold_finite_2006,arnold_finite_2010} has also initiated a change of paradigm.  The main idea behind is that geometrical and topological properties of differential operators are key points to understand how their discrete counterpart can be derived. The finite element exterior calculus (FEEC) is the result of this work and aims at studying approximations of PDEs that arise from Hilbert complexes. Let  $W_1, W_2$ be Hilbert spaces along with a differential map $d:W_1\rightarrow W_2$. The fundamental data of FEEC is then given by the polynomial subspaces $W_1^h$ and $W_2^h$ determined by projection maps $\pi_1$ and $\pi_2$ such that the following diagram commutes:
\begin{center}
    \begin{tikzcd}[cramped, column sep=large, row sep=huge]
      W_1 \arrow["d",r] \arrow["\pi_1"',d]  & W_2 \arrow["\pi_2",d] \\
      W_1^h \arrow["d",r]  & W_2^h
    \end{tikzcd}
\end{center}
The discretization can be again summarized by the triple $(W,d,\Omega(W))$ where $W$ is a polynomial algebra, $d$ a derivation map generating the exterior algebra $\Omega(W)$ with coefficients in $W$.\par
One can also mention of Christiansen et al. \cite{christiansen_construction_2008} on compatible differential forms on simplicial complexes.\par
Geometric integration and more generally structure preserving methods have applied this change of paradigm too \cite{arnold_compatible_2006,olver_geometric_2001,marsden_discrete_2001,hairer_geometric_2006,christiansen_topics_2011}. Symmetries and conservation laws of discrete operators parallel their continuous counterparts \cite{hydon_variational_2004,wan_multiplier_2016}. It has been shown that long-term stability can be obtained as a by-product \cite{wan_arbitrarily_2018}.  Finally, for an application of Lie groups to construct invariant discretization schemes, one can refer to \cite{bihlo_invariant_2013}.

Overall, in the geometric discretization framework, the realization is that classical analysis of consistency and stability is no longer the main criteria to look for in a discretization. In that context, consistency and stability are a consequence of preserving geometrical properties.

\paragraph*{Present work} The main question that we would like to address in this work is the existence of a unifying framework to geometric discretizations. This question can be divided into three subsidiary questions. \\
\textbf{The space:} the existence of a sequence of approximating spaces, with topological structures and metric specified at an early stage, that converges — in a suitable sense — to a manifold. \\
\textbf{The algebra:} tied to the question of space is the question of the algebra of “functions” and local coordinates. One needs to identify an associative algebra playing the role of the algebra of continuous functions over a space that do not necessarily possess a manifold structure.  It is a well established fact from the theory of Banach algebras \cite{elliott_simple_1974} that $C^*$-algebras can be realized as the set of continuous sections over some topological space. In a very intuitive description, an element of a $C^*$-algebras can be thought as (noncommutative) functions over a space called the spectrum of a $C^*$-algebra \cite{dixmier_c-algebras_1982,blackadar_operator_2006}. Hence, if one identifies the points of this spectrum one-to-one with the usual points of a topological space $X$, then a $C^*$-algebra appears as a good candidate for the set of continuous functions over $X$. Indeed, their normed space structure is a powerful tool to study boundedness and convergence of its elements. Thus, in the same fashion as for the space itself, can one construct a nested sequence of algebras such that the limit is essentially the space of continuous functions over the original manifold ? \\
\textbf{The geometry:} once the questions of space and algebra are addressed, it remains to define (if it exists and is it unique ?) a differential calculus — understood from an algebraic/geometric point of view in opposition to the usual analytic perspective — on such space. What does such a differential structure on an infinitesimal space look like ? One can already notice that it will irremediably differ from its continuous counterpart since functions and forms do not commute anymore: 
\vspace{0.5cm}
\begin{equation}
gdf\neq dfg.
\end{equation}
Moreover, the differential calculus is intimately tied to connections and distances between points parametrized by $h$. This fact is reminiscent of the continuous theory, where the line element $ds$ — one can think of an infinitesimal displacement vector in a metric space — on a $n$-dimensional Riemannian manifold is a function of the metric tensor. Moreover, it is a well established fact in spin geometry \cite[pp. 552-557]{connes_noncommutative_1994} that this metric information can be summarized in a single operator $\slashed{D}$ called the Dirac operator \cite[pp. 406-407]{rudolph_differential_2017} such that :
\begin{equation}
ds =\slashed{D}^{-1}.
\end{equation}
Therefore, topology, metric, and differentiation can be deduced — in principle — from the data of the Dirac operator. Hence, one have a dual description of space : one purely topological given by an open cover and one purely algebraic given by the Dirac operator.
\paragraph*{Objectives} The main objective of this work is to derive \textit{ab initio} finite difference calculus using the language of noncommutative geometry. This leads us to define tools from differential geometry such as differential maps along with their differential complex, affine connections and a Laplace operator. It also allows us to study spectral convergence with respect to a parameter $h$. Indeed, the natural setting of $C^*$-algebras, through their representations into operator algebras, allows us to use the machinery of functional calculus. This main objective can be divided into three sub-objectives. First, we aim at establishing a proper notion of \textit{discrete space} $X$, starting uniquely from the knowledge of a manifold $M$ along with its algebra of functions. Secondly, we want to exhibit the algebra of continuous sections $\Gamma(X)$ over $X$. Following Gelfand-Naimark's theorem, this should be a $C^*$-algebra $A$. Thirdly, we define a so-called \textit{Dirac operator} $D$ governing the differential geometry over the space $X$. Once such operator is defined, it provides an exterior algebra $\Omega(A)$ and some usual machinery from differential geometry.\\
In this work, we are able to give an intrinsic description of finite difference calculus in terms of noncommutative geometry and its quantized calculus. We recover some usual tools of differential geometry, such as an exterior derivative. Higher-order approximations are restated in terms of $\mathbb{Z}_2$-graded traces induced by positive operators. We also define and establish convergence of differential operators on infinitesimal spaces to their continuous counterpart. We further prove a generalized result on direct limits of $C^*$-algebras over posets. This extends the result of Bimonte et al. \cite{bimonte_noncommutative_1996} proven in the special case of noncommutative lattice.
Therefore, this work opens the door to a general framework to study approximation theory of PDEs.
\paragraph{Technical formalism} In this paper we consider the formalism of Noncommutative Differential Geometry (NDG). NDG has been introduced by Connes in a series of papers \cite{connes_non-commutative_1985} compiled in the red book \cite{connes_noncommutative_1994} — and later extensively developed by Connes and his collaborators \cite{connes_noncommutative_2007,connes_noncommutative_2019}. This branch of mathematics is concerned with a geometric approach to noncommutative algebras \cite{khalkhali_very_2004, varilly_introduction_2006,sergeev_intorduction_2016,sitarz_lectures_2008,sitarz_friendly_2013}. In Connes' work, a noncommutative space is — heuristically speaking — the dual space of a $C^*$-algebra by analogy to the Gelfand theory where commutative $C^*$-algebras are dual objects to locally compact Hausdorff spaces. In fact, the notion of space becomes secondary and is replaced by the notion of a spectral triple $(A,H,D)$ — where $A$ is a $C^*$-algebra, $H$ is a Hilbert space on which $A$ is realized as an algebra of bounded linear operators, $D$ is a Dirac operator. A new type of differential calculus using functional analysis is then derived; it is now referred to as quantized calculus. We also mention another type of noncommutative differential geometry over matrix algebras developed by Dubois-Violette et al. \cite{duboisviolette_noncommutative_1990} 
and exposed in more details in \cite{madore_introduction_1999,masson_gauge_2012,masson_geometrie_1996}.\par
The idea of approximating a bounded region of space-time with finite topological spaces as been pushed by Sorkin \cite{sorkin_finitary_1991}.
Important examples of noncommutative spaces are provided by noncommutative lattices, which are a particular case of posets. This topic has been thoroughly studied by Bimonte et al. \cite{bimonte_lattices_1996,bimonte_noncommutative_1996} — summarized in Landi's book \cite{landi_introduction_1997} —
and techniques from noncommutative geometry have been used to construct models of gauge theory on these noncommutative lattices in Balachandran et al. \cite{balachandran_noncommutative_1996,balachandran_lattice_1998,balachandran_classical_2001}. It is also worth mentioning another approach to discrete noncommutative spaces and their differential calculus in the work of Dimakis et al. \cite{dimakis_discrete_1994}.\\

In this paper, we start by reviewing some technical preliminaries in Section \ref{Prelem}, then we construct $C^*$-algebra over posets in Section \ref{Calg}. We then construct the differential structure in Section \ref{Gem} and conclude with a discussion on the convergence results for the $d$-lattice in $\mathbb{R}^d$ and for the torus $\mathbb{T}^d$.       

\section{Preliminaries}
\label{Prelem}
Unless stated otherwise, we will consider $M$ to be a smooth compact connected manifold $M$ of dimension $d$.
\subsection{Triangulation and posets}
\begin{figure}[ht]
\centering
\adjustbox{trim=3cm 10cm 3cm 3cm}{
\resizebox{120mm}{!}{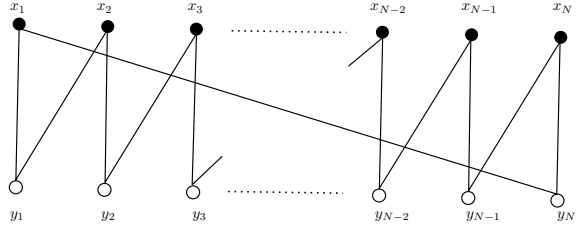}}
\caption{Poset associated to a triangulation of $S^1$.}
 \label{fig:Ncircle}
\end{figure}

Let $K$ be an \textit{abstract simplicial complex} with elements $\sigma$ and $|K|$ its geometric realization. The dimension of a simplex $\sigma \in |K|$, denoted $\dim(\sigma)$, is the dimension of the smallest affine space containing $\sigma$. The set $K$ can be written as a union of subsets $K(n)$, where $\sigma^n \in K(n)$ is a simplex of dimension $n$. The subset $K(0)$, also denoted $V$, is the set of vertices; the subset $K(1)$, also denoted $E$, is the set of edges.\\
A manifold $M$ admits a \textit{triangulation} $\mathcal{T}(K)$ if there exists a simplicial complex $K$ and homeomorphism $\varphi :|K|\rightarrow M$ between $M$ and the geometric realization $|K|$. We recall the following theorem due to Whitney on the existence of a triangulation.
\begin{mth}[ {\cite[pp.124-135]{whitney_geometric_1957}} ] Every $k$-smooth manifold $M$ admits a triangulation, for $k\geq 1$.
\label{thm2}
\end{mth}
To every simplicial complex $K$, one can associate a partially ordered set (poset) $P(K)$ which is defined to be the poset of nonempty faces ordered by inclusion. We will denote by $\leq$ the partial order on $P(K)$. The preorder $\leq$ induces a topology $P(K)$ called the \textit{Alexandrov topology} and generated by the bases of open sets $\mathcal{B}:=\left\lbrace U_x:=\left\lbrace y\in X : y\leq x \right\rbrace : x\in X \right\rbrace $. For instance, Figure \ref{fig:Ncircle} show the poset associate to a triangulation of the circle $S^1$. \par
Reciprocally, to every poset $X$, one can associate an abstract simplicial complex $K(X)$, where the simplices are nonempty chains in $X$.\\
A map $f:X\rightarrow Y$ between posets is continuous if and only if it is order preserving with respect to the orders associated with the order topologies on $X$ and $Y$. The map $f$ induces a simplicial map $K(f):K(X)\rightarrow K(Y)$; vice-versa to every simplicial map $f:K\rightarrow L$, one can associate a continuous map $P(f):P(K)\rightarrow P(L)$ between posets.\\
One can reverse the order $\leq$ on a poset $X$ and define the space $X^{op}$. These spaces have the same underlying set. Open sets in $X$ correspond to closed sets in $X^{op}$ and vice-versa. Moreover, a continuous map $f:X\rightarrow Y$ induces a continuous map $f^{op}: X^{op}\rightarrow Y^{op}$ and vice-versa.\\
Finally, the space $K(P(K))$ is called \textit{the barycentric subdivision} of the simplicial complex $K$ and is denoted $K'$. In addition, $K'$ is a simplicial complex and there exists a continuous embedding $i:K'\rightarrow K$. It identifies $K'$ as a subspace of $K$. Furthermore, the map $i$ also induces a continuous embedding on the posets:
\begin{equation*}
P(i):P(K')\rightarrow P(K),
\end{equation*}
where the elements of $P(K')$ are nonempty chains of $P(K)$. If the complex $K$ is in a metric space, then one can define the diameter $\text{diam}(\sigma)$ of a simplex $\sigma$; the largest of these is the mesh of $K$. We can then inductively form the $n$-th barycentric subdivision $K_n=(K_{n-1})'$;  the sequence $(K_n)$ can be constructed such that $mesh(K_n)\rightarrow 0$. We will denote by $h_n$ the mesh length of $K_n$.\\
In the rest of this work, we will consider the space $X_n=P(K_n)^{op}$ where the elements are the simplices of $K_n$ and the ordering is by reversed inclusion. The poset $X_n$ is equipped with the Alexandroff topology induced by the inclusion order. Starting from a triangulation $\mathcal{T}(K)$ of $M$ and a homeomorphism
\begin{equation*}
\varphi : |K|\rightarrow M,
\end{equation*}
 we construct a sequence of posets $(X_n)$ associated to the successive barycentric subdivisions $(K_n)$ of $K$. The maps $\phi_{n,m}:X_m\rightarrow X_n$ for $m\geq n$ sending an element from $X_m$ to its carrier in $K_n$ form a sequence $\left\lbrace X_n,\mathbb{N},\phi_{n,m} \right\rbrace $:
\begin{center}
\begin{tikzcd}
X_0 & X_1 \arrow["\phi_{12}"',l] & X_2 \arrow["\phi_{23}"',l] & X_3 \arrow["\phi_{34}"',l] & \cdots \arrow["\phi_{45}"',l]
\end{tikzcd}
\end{center}
\begin{figure}[ht]
\centering
\adjustbox{trim=3cm 11cm 3cm 0cm}{
\resizebox{120mm}{!}{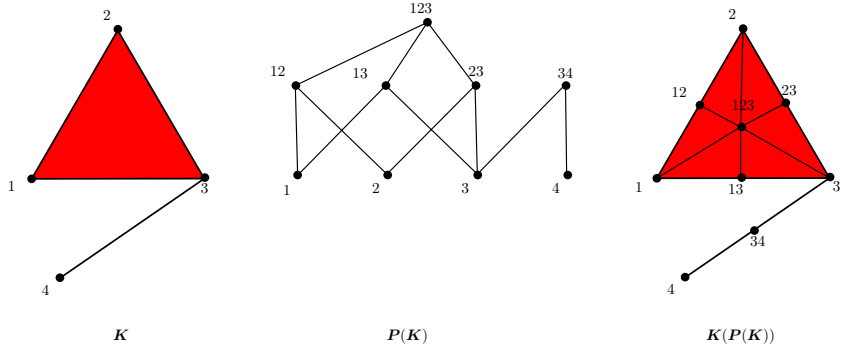}}
\caption{Simplicial complex, poset and barycentric subdivision.}
 \label{fig:poset}
\end{figure}
\subsection{The inverse limit construction}
\label{inverse}
We have the system $\left\lbrace X_n,\mathbb{N},\phi_{n,m} \right\rbrace $ where the maps $\phi_{n,m}$ satisfy by construction the coherence properties, for $\leq n \leq m$:
\begin{equation}
\phi_{l,n}\circ \phi_{n,m}=\phi_{l,m},  \quad \phi_{n,n}=id.
\end{equation}
Therefore, the system $\left\lbrace X_n,\mathbb{N},\phi_{n,m} \right\rbrace $ defines an inverse system of topological spaces. We define its inverse limit 
\begin{equation}
X_{\infty}:=\lim_{\leftarrow}X_i 
\end{equation}
which is a subset of the product space ${\prod}_{i\in I}X_i$ and we topologize it with the subspace topology. An element $x\in X_\infty$ is then a coherent sequence:
\begin{equation}
x=(x_1,x_2,\cdots,x_i, \cdots x_j, \cdots )\in {\prod}_{i\in I}X_i, \quad x_i=\phi_{i,j}(x_j) \  \forall i\leq j .
\end{equation}
Equivalently, recalling the definition of $X_n$ from a simplicial complex $K_n$, one can see an element of $X_\infty$ as a coherent sequence of nested simplices. The inverse limit $X_\infty$ also comes equipped with natural projection maps $\phi_i:X_\infty \rightarrow X_i$ which pick out the $i$-th coordinate for every $i\in \mathbb{N}$.\\
The space $X_\infty$ is a poset; the partial order on the sets $X_n$ give a partial order $\leq$ on the set $X_\infty$, where $y\leq x$ provided that $y_n\leq x_n$ for every $n\in \mathbb{N}$. Moreover, using the homeomorphism between $M$ and $|K|$, we see that there is a natural map $p_n:M\rightarrow X_n$ for each $n$, since every point in $K$ is contained in the interior of exactly one face of the $n$-th barycentric subdivision of $K$. We have the following commuting diagram:
\begin{center}
\begin{tikzcd}
{} & M \arrow["p_{n-1}"',ld] \arrow["p_{n}",rd] \\
X_{n-1} & {} & X_{n} \arrow["\phi_{n-1,n}",ll] 
\end{tikzcd}
\end{center}
In addition, using the correspondence between points in $X_n$ and faces of simplices in $K_n$, we can denote the simplex corresponding to $x_n\in X_n$ by $\sigma_n(x)$. We then immediately have that for every $n\geq 0$:
\begin{equation*}
p_n^{-1}(U_x)=\text{st}(\sigma_n(x)),
\end{equation*}
where $\text{st}$ is the open star map. This implies that the maps $p_n$ are continuous. We can then define a continuous map 
\begin{equation}
p:M\rightarrow X_\infty, \quad p(a) = (p_0(a), p_1(a), \cdots).
\end{equation}
The next claim allows us to create a map from $X_\infty$ to $M$ which acts as an inverse to $p$.
\begin{mlem}
Given $x=(x_0,x_1,\cdots)\in X_\infty$, pick $a_n\in p_n^{-1}(x_n)$ for each $n\geq 0$. Then the sequence $(a_n)$ converges to $a\in M$ and the map 
\begin{equation*}
G:X_\infty \rightarrow M, \quad x \mapsto a_x
\end{equation*}
is well-defined and continuous.
\label{Gcont}
\end{mlem}
\begin{proof}
The points $a_n \in K_n$ lie in nested simplices of increasingly fine barycentric subdivisions of $K$. Any sequence obtained this way converges to the same point since they are obtained by intersection of nested closed sets with vanishing diameters. The proof of continuity of $G$ can be found in \cite[Prop. 2.4.16]{thibault_homotopy_2013}.
\end{proof}
\begin{mlem} Let $x \in X_\infty$ such that $G(x)=a_x$, then $p(a_x)\geq x$.
\label{ineq}
\end{mlem}
\begin{proof}
Recall that the order in $X_\infty$ is given by: $x\leq y$ in $X_\infty$ if and only if $x_n\leq y_n$ in $X_n$ for every $n$.\\
Now suppose that $p(a_x)\geq x$ is not true, then there exists $n$ such that $p(a_n)\geq x_n$ is not true. This means that $p(a_n)$ is not contained in the simplex corresponding to $x_n\in X_n$. Thus, it contradicts the fact that $a_n\in p_n^{-1}(x_n)$.
\end{proof}
\begin{mlem}
The set $p(M)$ is precisely the subspace $\mathfrak{M}$ of all maximal elements in $X_\infty$.
\end{mlem}
\begin{proof}
Let $y$ be a maximal element in $X_\infty$. Then by Lemma \ref{ineq}, $p(a_y)\geq y$ and therefore $p(a_y)=y$. Conversely, if there exists $a\in M$ and $y\in X_\infty$ such that $y\geq p(a)$, then by definition, $y_n\geq p_n(a)$ for every $n$. Now, let $G(p_n(a))=a_n$ and $G(y)=y_n$ for every $n$. Because $y_n\geq p_n(a)$, we have that $y_n\in p^{-1}_n(p_n(a))$ for every $n$. Hence, the sequences $(a_n)$ and $(y_n)$ have the same limit $a_y=a$. Thus, $p(a_y)=p(a)$ and $p(a)\geq y$ again by Lemma \ref{ineq}. We conclude that $p(a)=y$. 
\end{proof}
\begin{mprop} The space $M$ is homeomorphic to the subspace $\mathfrak{M}$ of all maximal points of the inverse limit of the system $\left\lbrace X_n,\mathbb{N},\phi_{n,m} \right\rbrace $.
\label{main0}
\end{mprop}
\begin{proof}
We need to prove that $G:p(M)\rightarrow M$ is a homeomorphism. By construction, we have that $G\circ p = id$, then, $G$ is a bijection. By Lemma \ref{Gcont}, $G$ is continuous. \\
Since $p(M)$ is equipped with the subspace topology, an open set $U$ pf $p(M)$ can be written as $U=V\cap p(M)$ where $V$ is an open set in $X_\infty$. Now $G(U)=p^{-1}(V)$, thus $G(U)$ is open. Hence, $G$ is a continuous and open bijective map and thus a homeomorphism.
\end{proof}
\subsection{C*-algebras and their spectra}
We conclude this section by introducing some of the fundamental concepts on $C^*$-algebras that will be useful in the rest of this paper; thorougher details can be found in the literature \cite{dixmier_c-algebras_1982,blackadar_operator_2006,murphy_c-algebras_1990}.

A $C^*$- algebra $A$ is a Banach algebra over $\mathbb{C}$ together with an involution $x\mapsto x^*$ such that:
\begin{equation}
\begin{array}{ccccc}
(xy)^*=y^*x^* &\text{and} & \|x^*x\|=\|x\|^2 & \text{for} & x,y\in A. 
\end{array} 
\end{equation}
The two archetypes of $C^*$-algebras are given by the space of continuous complex-valued functions that vanish at infinity $\left( C_b(X),\|\cdot\|_{\infty}\right) $ over a locally compact Hausdorff space $X$ — in the commutative setting — and the space of bounded operators $\left( B(H),\|\cdot\|_{\text{op}}\right) $ over a Hilbert space $H$ — in the noncommutative case.\\
A central tool in the study of $C^*$-algebras is through their representations.
\begin{mdef}[Representation] Let $A$ be a $*$-algebra. A representation of $A$ is a pair $(\pi,\mathcal{H})$ where $\mathcal{H}$ is a Hilbert space and $\pi:A\rightarrow \mathcal{B}(\mathcal{H})$ is a $*$-homomorphism. We also say that $\pi$ is a representation of $A$ on $\mathcal{H}$.
\end{mdef} 
Another crucial tool to study $C^*$-algebras, and related to their representations, is the primitive spectrum.
\begin{mdef}[Primitive spectrum]
The \textit{primitive spectrum} $Prim(A)$ is the space of kernels of irreducible $^*$-representations equipped with the \textit{hull-kernel (Jacobson) topology}.
\end{mdef}
The primitive spectrum becomes central to describe the internal algebraic structure of $A$. It can be turned into a topological space using the \textit{hull-kernel (Jacobson) topology}. Let $W\in 2^{Prim(A)}$ an element of the power set, then the closure operator is given by
\begin{equation*}
Cl(W):= \left\lbrace I\in Prim(A) : \bigcap\ker(\pi)\subseteq I\right\rbrace.
\end{equation*}
A related and equally important notion, is the \textit{spectrum} $Spec(A)$ of a $C^*$-algebra i.e. the set of non-zero unitary equivalence classes of irreducible $^{*}$-representations. There is an immediate surjection map 
\begin{equation}
Spec(A)\rightarrow Prim(A), \quad (\mathcal{H},\pi)\mapsto \ker \pi,
\label{specprim}
\end{equation}
which endows $Spec(A)$ with the pull-back of the Jacobson topology.
\begin{mrmk} When the primitive spectrum $Prim(A)$ is a $T_0$-space, then the map \eqref{specprim} is a homeomorphism. This will always be the case in this work, therefore we will indifferently refer to the primitive spectrum or to the spectrum.
\end{mrmk}
In the commutative case, the spectrum of $A$ plays the role of a \textit{space}. Indeed, any element $a\in A$ can be interpreted as a function over the space of characters through the \textit{Gelfand map}:
\begin{equation}
a\ni A\mapsto (\chi \mapsto\hat{a}(\chi) ) \ (\chi \in Spec(A)).
\end{equation}
If we let $X=Spec(A)$, then the Gelfand transform is an isomorphism of $A$ onto the $C^*$-algebra $C(X)$ of continuous complex functions over $X$.
\section{C*-algebras over a triangulation}
\label{Calg}
In this section, we show how to associate a $C^*$-algebra $A_n$ to the space $X_n$ defined in the previous section. The construction follows the works of Behncke and Leptin \cite{behncke_c-algebren_1972, behncke_c-algebras_1972, behncke_c-algebras_1973, behncke_class_1973}. In order to give a more comprehensive presentation, we states the procedure as a sequence of axioms in the subsection \ref{BehLep}. For more details, we refer to \cite{ercolessi_noncommutative_1996}.\\
In the rest of this work, $A$ will designate a $C^*$-algebra (eventually infinite dimensional) and $H$ a representation of $A$. The letters $\mathfrak{A}$ and $\mathfrak{H}$ will be used in the commutative case.
\subsection{C*-algebras over a topological space}
We let $X$ be a topological space. A $C^*$-algebra over $X$ is a pair $(A,\psi)$ consisting of a $C^*$-algebra $A$ and a continuous surjection
\begin{equation*}
\psi:Prim(A)\rightarrow X.
\end{equation*}
Let $\mathcal{O}_X$ be the set of open subsets of $X$, partially ordered by inclusion. For a $C^*$-algebra $A$, we let $I(A)$ be the set of all closed $^*$-ideals in $A$ partially ordered by inclusion. There is an isomorphism (see \cite{meyer_c-algebras_2008}) between $I(A)$ and the set of open subsets $\mathcal{O}_{Prim(A)}$ in $Prim(A)$. We will always identify $\mathcal{O}_{Prim(A)}$ and $I(A)$ through the isomorphism:
\begin{equation}
\mathcal{O}_{Prim(A)}\simeq I(A) \qquad U\mapsto \bigcap_{\pi\in Prim(A)\backslash U}\pi.
\label{openideal}
\end{equation}
Then for $(A,\psi)$ a $C^*$-algebra over $X$, we get a map
\begin{equation*}
\psi^*:\mathcal{O}_X\rightarrow \mathcal{O}_{Prim(A)}\simeq I(A) \qquad U\mapsto \lbrace \pi\in Prim(A) | \psi(\pi)\in U\rbrace \simeq A(U).
\end{equation*}
We will denote by $A(U)\in I(A)$ the ideal associated to the open subset $U$. We can now identify the open sets in $X$ with closed $^*$-ideals of $A$, and points in $X$ with irreducible representations of $A$.
\subsubsection{The Behncke-Leptin construction}
\label{BehLep}
The Behnck-Leptin construction allows us to associate a $C^*$-algebra $(A,\psi)$ over a partially ordered space $X$ such that $\psi=id$ is the identity map. Hence, the spaces $Prim(A)$ and $X$ can be identified.\\

The axioms of the Behncke-Leptin construction go as follows:
\begin{itemize}
\item[1)] Associate a separable Hilbert space $H(X)$ to the space $X$ and attach to every point $x\in X$ a subspace $H(x)\subseteq H(X)$ that decomposes into:
\begin{equation}
H(x)=H^-(x)\otimes H^+(x).
\end{equation}
where $H^-(x)\simeq \ell^2(\mathbb{Z)}$.
\item[2)] Let $\mathfrak{M}$ be the set of maximal points in $X$. Then for every $x\in \mathfrak{M}$, one has 
\begin{equation}
H(x)=H^-(x)\otimes \mathbb{C}\simeq H^-(x).
\end{equation}
\item[2')] If $\mathfrak{m}$ is the set of minimal points in $X$, then for every $x\in \mathfrak{m}$, one has 
\begin{equation}
H(x)=\mathbb{C}\otimes H^+(x)\simeq H^+(x).
\end{equation}
\item[3)] Associate to each point $x\in X$ an operator algebra $A(x)$ acting on $H(x)$ (extended by zero to the whole space $H(X)$) such that
\begin{equation}
A(x)=1_{H^-(x)}\otimes \mathcal{K}(H^+(x)).
\label{subalgebra}
\end{equation}
where $\mathcal{K}(H^+(x))$ is the set of compact operators over $H^+(x)$.
\item[4)] Build the $C^*$-algebra $A(X)$ associated to the space $X$ as the algebra generated by the subalgebras $A(x)$ when $x$ run over $X$:
\begin{equation}
A(X)=\bigoplus_{x\in X}A(x) \hspace{0.3cm} \text{acting on} \hspace{0.3cm} H(X)=\bigoplus_{x\in X}H(x).
\end{equation}
\end{itemize}
As mentioned already, using the isomorphism 
\begin{equation}
\psi:X\rightarrow Spec(A), \quad \psi(x)=\pi_x
\label{iso}
\end{equation}
one can identify a point $x\in X$ with an irreducible representation $(\mathcal{H}_x,\pi_x)$:
\begin{equation}
\pi_x:A(X)\rightarrow B(\mathcal{H}_x), \quad a\mapsto \pi_x(a).
\end{equation}
The irreducible representation $\mathcal{H}_x\subset H(x)$ is obtained as a subspace of $H(x)$. We define the following total space:
\begin{equation}
\mathcal{H}_X=\bigoplus_{x\in X}\mathcal{H}_x.
\label{representation}
\end{equation}
An element $a\in A$ then uniquely defines a \textit{map} on $X$:
\begin{equation}
\hat{a}:X\rightarrow A, \quad \hat{a}(x):=\pi_x(a)=\sum_{i\in I_x}\lambda_i(x)1\otimes k_i(x)
\label{map}
\end{equation}
where $\lambda_i(x)\in \mathbb{C}$ and $k_i(x)$ is a compact operator. In particular, if we identify the Hilbert space $H^{-}(x)$ with $\ell^2(\mathbb{Z})$, then we see that $\lambda(x) 1$ is nothing else than a multiplication operator: 
\begin{equation}
T_{\lambda(x)} (u) = \lambda(x) \cdot u.
\end{equation}
for $u\in \ell^2(\mathbb{Z})$. This leads us to the fifth axiom.
\begin{itemize}
\item[5)] For every $x\in \mathfrak{M}$, the representation $(\mathcal{H}_x,\pi_x)$ is one-dimensional:
\begin{equation}
\pi_x:A(X)\rightarrow \mathbb{C}, \quad a\mapsto \pi_x(a)=\lambda(x).
\end{equation}
\end{itemize}
\begin{mex} Let $\sigma$ be a $2$-simplex and consider $X$ to be the poset associated to $\sigma$ with the opposite order. Then, Figure \ref{fig:algebra} shows a generic element $a_x\in A(x)$ for every vertex $x$ of $X$. The full algebra $A(X)$ is obtained as a direct sum of the $A(x)$.
\begin{figure}[ht!]
\centering
\adjustbox{trim=3cm 7cm 3.5cm 1cm}{
\resizebox{100mm}{!}{
\begingroup%
  \makeatletter%
  \providecommand\color[2][]{%
    \errmessage{(Inkscape) Color is used for the text in Inkscape, but the package 'color.sty' is not loaded}%
    \renewcommand\color[2][]{}%
  }%
  \providecommand\transparent[1]{%
    \errmessage{(Inkscape) Transparency is used (non-zero) for the text in Inkscape, but the package 'transparent.sty' is not loaded}%
    \renewcommand\transparent[1]{}%
  }%
  \providecommand\rotatebox[2]{#2}%
  \newcommand*\fsize{\dimexpr\f@size pt\relax}%
  \newcommand*\lineheight[1]{\fontsize{\fsize}{#1\fsize}\selectfont}%
  \ifx\svgwidth\undefined%
    \setlength{\unitlength}{595.27559055bp}%
    \ifx\svgscale\undefined%
      \relax%
    \else%
      \setlength{\unitlength}{\unitlength * \real{\svgscale}}%
    \fi%
  \else%
    \setlength{\unitlength}{\svgwidth}%
  \fi%
  \global\let\svgwidth\undefined%
  \global\let\svgscale\undefined%
  \makeatother%
  \begin{picture}(1,1.41428571)%
    \lineheight{1}%
    \setlength\tabcolsep{0pt}%
    \put(0,0){\includegraphics[width=\unitlength,page=1]{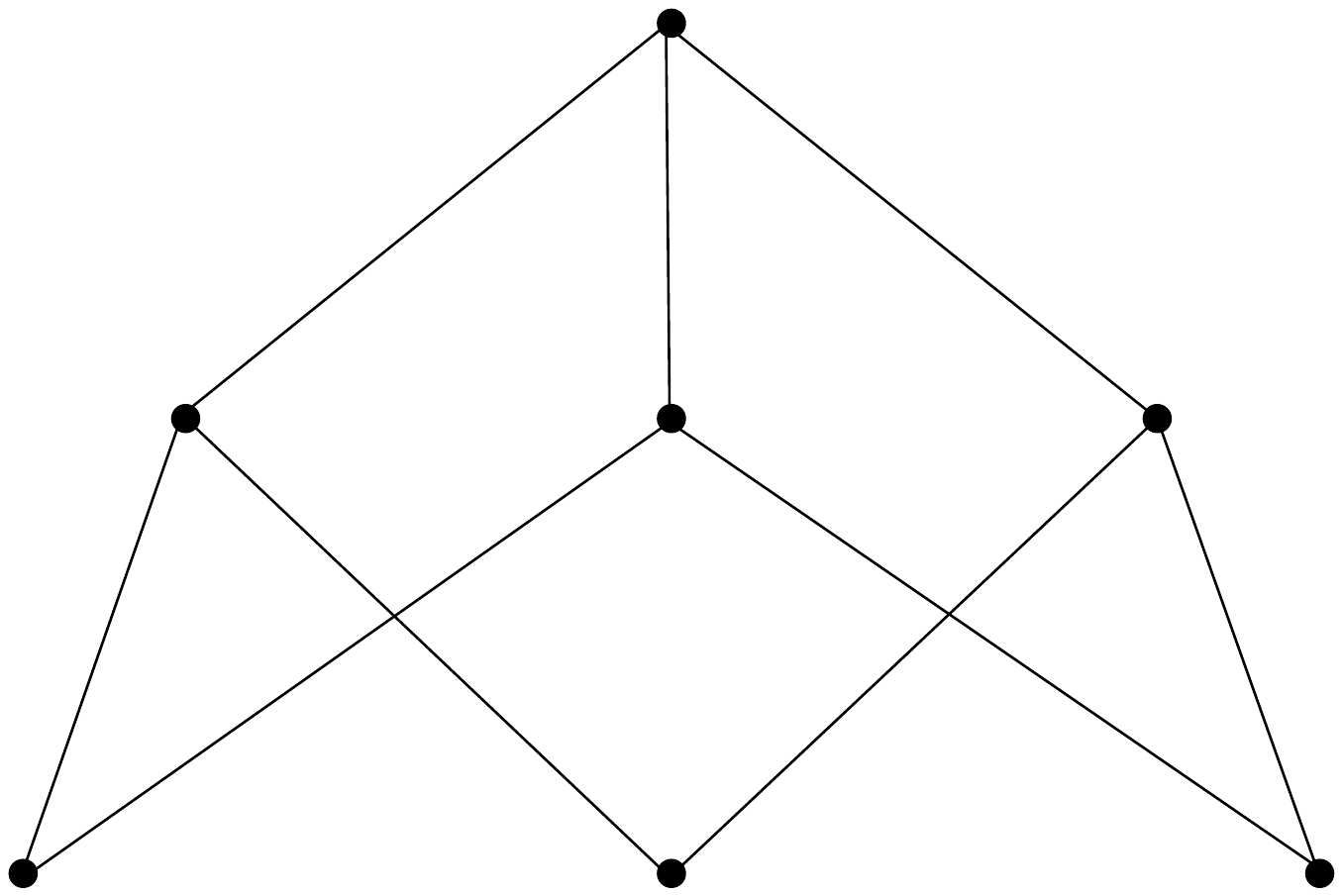}}%
    \put(0.46400586,1.25197193){\color[rgb]{0,0,0}\makebox(0,0)[lt]{\lineheight{1.25}\smash{\begin{tabular}[t]{l}\Large$\lambda(x)1$ \end{tabular}}}}%
    \put(0.09,1.06){\color[rgb]{0,0,0}\makebox(0,0)[lt]{\lineheight{1.25}\smash{\begin{tabular}[t]{l}\Large$\mu(y_1)1\otimes k(y_1)$\end{tabular}}}}%
    \put(0.46753438,1.06){\color[rgb]{0,0,0}\makebox(0,0)[lt]{\lineheight{1.25}\smash{\begin{tabular}[t]{l}\Large$\mu(y_2)1\otimes k(y_2)$\end{tabular}}}}%
    \put(0.7,1.06){\color[rgb]{0,0,0}\makebox(0,0)[lt]{\lineheight{1.25}\smash{\begin{tabular}[t]{l}\Large$\mu(y_3)1\otimes k(y_3)$\end{tabular}}}}%
    \put(0.08,0.8){\color[rgb]{0,0,0}\makebox(0,0)[lt]{\lineheight{1.25}\smash{\begin{tabular}[t]{l}\Large$k(z_1)$\end{tabular}}}}%
    \put(0.47988434,0.8){\color[rgb]{0,0,0}\makebox(0,0)[lt]{\lineheight{1.25}\smash{\begin{tabular}[t]{l}\Large$k(z_2)$\end{tabular}}}}%
    \put(0.79,0.8){\color[rgb]{0,0,0}\makebox(0,0)[lt]{\lineheight{1.25}\smash{\begin{tabular}[t]{l}\Large$k(z_3)$\end{tabular}}}}%
  \end{picture}%
\endgroup%
}}
\caption{$C^*$-algebra associated to a poset.}
 \label{fig:algebra}
\end{figure}
\end{mex}
\subsubsection{Commutative subalgebras}
Let $(A,id)$ be the $C^*$-algebra associated to a finite connected poset $X$ through the Behncke-Leptin construction. Among the subalgebras of $A$, those of particular interest are commutative ones. The centre of $A$ will be denoted by $Z(A)$. We know by construction that $A$ is generated by the algebras 
\begin{equation*}
A(x)=1_{H^-(x)}\otimes \mathcal{K}(H^+(x))
\end{equation*}
for $x$ running $X$. Moreover, we recall that the algebra of compact operators $\mathcal{K}(H)$ over an infinite dimensional Hilbert space $H$ has a trivial centre. We deduce that, for a given $x\in X$, $A(x)$ has a trivial centre.
\begin{mprop} The centre $Z(A)$, also denoted by $\mathfrak{A}$, of $A$ is given by
\begin{equation}
\mathfrak{A}=\oplus_{x\in \mathfrak{M}}1_{H(x)}.
\label{commutative}
\end{equation}
and is generated by the projectors on $H(x)$ where $x\in \mathfrak{M}$ is a maximal point.
\end{mprop}
\begin{proof} This is a direct consequence of the fact that the centre of $K(H)$ is trivial and the definition of the generating subalgebras $A(x)$ in the Behncke-Leptin construction.
\end{proof}
\subsection{C*-algebra over a simplicial complex}
We go back now to a simplicial complex $K$ and its associated poset $P(K)^{op}$ that we will denote $X$ (seen as a topological space). Using the Behncke-Leptin construction, we can associate a $C^*$-algebra $(A(X),id)$ over $X$ such that $
Prim(A)$ is  identified with $X$.\\
Now, let $K$ and $K'$ be simplicial complex such that $K'$ is a barycentric subdivision of $K$. We denote by $X$ and $X'$ the associated posets. We then have a continuous surjection :
\begin{equation*}
\phi : X'\rightarrow X.
\end{equation*} 
Consider in addition that $(A(X),id)$, respectively $(A(X'),id')$, is a $C^*$-algebra over $X$, respectively $X'$.  We would like to show that for the given map $\phi$, there exists a pullback map $\phi^*$ such that the following diagram commutes:
\begin{center}
\begin{tikzcd}[cramped, column sep=large, row sep=huge]
A(X) \arrow["\phi^*",r]\arrow["id"',d] &  A(X') \arrow["id'",d]\\
X & X' \arrow["\phi",l]
\end{tikzcd}
\end{center}
i.e. such that the following proposition is satisfied:
\begin{equation}
\pi_x(a) = \pi_{y}(\phi^*(a)), \quad \forall x\in X, \ \forall \text{$y\in \phi^{-1}(x)$ : $\dim(\sigma'_y)=\dim(\sigma_x)$}.
\label{pullback}
\end{equation}
Here, $\sigma_x$ and $\sigma'_y$ are the simplex associated to $y$ and $x$ in the identification of $X$ and $X'$ with $K$ and $K'$. We are also using the isomorphism \eqref{iso} to identify a point $x\in X$ with an irreducible representation $(H_x,\pi_x)\in Spec(A)$; then $\pi_x(a)$ is an operator acting on $H_y$ and $\pi_{y}(\phi^*(a))$ an operator on $H'_y$. We are assuming here that $H_x$ and $H'_y$ can be identified as Hilbert spaces; the identification is constructed in Equation \eqref{hilbert}.
\label{pullback2}
\begin{mprop}
A continuous surjection $\phi:X'\rightarrow X$ between posets induces a unital $^{*}$-homomorphism $\phi^*:A(X)\rightarrow A(X')$ satisfying \eqref{pullback}.
\end{mprop}
\begin{proof}
We recall that the algebra $A(X)$ is generated by the subalgebras $A(x)$ defined by \eqref{subalgebra} for $x$ running in $X$. Then, it is enough to define $\phi^*$ on the algebras $A(x)$ and extend the map by linearity.\\

Therefore, if we start with the following decomposition:
\begin{equation}
A(X_i)=\oplus_{x\in X_i}A(x), \quad a=\sum_{x\in X}a_x, \quad \text{i=1,2}
\end{equation}
with $X_1=X$ and $X_2=X'$, we define $\phi^*$ such that:
\begin{equation}
\phi^*(a)=\sum_{y\in X'}a_y,
\end{equation}
where 
\begin{equation}
a_y = \left\lbrace 
\begin{array}{cc}
a_\phi(y) & \text{if $\dim(\sigma_y)=\dim(\sigma_{\phi(y)})$,} \\ 
0 & \text{otherwise}.
\end{array} 
\right. 
\end{equation}
Thus, if we let $x\in X$ and consider the set:
\begin{equation}
\Phi^{-1}(x)=\lbrace y \in \phi^{-1}(x) : \dim(\sigma'_y)=\dim(\sigma_x) \rbrace.
\end{equation}
then, we have defined $\phi^*$ such that it satisfies \ref{pullback} i.e. for any $a\in A(X)$:
\begin{equation}
\pi_{y}(\phi^*(a)) = \pi_x(a), \quad \forall y\in \Phi^{-1}(x).
\end{equation}
Furthermore, $\phi^*$ is a $^*$-homomorphism by construction. In addition, the identity element on $A(X)$ is given by 
\begin{equation}
1_{A(X)}= \sum_{x\in \mathfrak{M}}1_{H(x)}
\end{equation}
and since $\phi(\mathfrak{M}')=\mathfrak{M}$, then $\phi^*(1_{A(X)})=1_{A(X')}$ i.e. $\phi^*$ is unital.
\end{proof}
\subsection{The direct limit construction}
We now recall the definition of a direct limit of $C^*$-algebras. Consider a direct sequence $(A_n,\psi_n)$ of separable $C^*$-algebras with *-homomorphism $\psi_n:A_n\rightarrow A_{n+1}$. The product $\prod_nA_n$ equipped with the pointwise addition, multiplication, scalar multiplication and involution is a $C^*$-algebra \cite{olai_milhoj_af-algebras_2018}. We denote by $A'$ the following set
\begin{equation}
A'=\left\lbrace a=(a_n)\in \prod_nA_n : \exists N\in \mathbb{N}, \ a_{n+1}=\psi_n(a_n)\  \forall n\geq N \right\rbrace.
\end{equation}
Since $(\psi_n)$ are contractions, then $(\|a_n\|)_n$ converges. One can then check that the map
\begin{equation}
p:A'\rightarrow \mathbb{R}^+, \ a\mapsto p(a):=\lim_{n\rightarrow \infty}\|a_n\|,
\end{equation}
is a $C^*$-seminorm on $A'$. The \textit{direct (inductive) limit} of the sequence $(A_n,\psi_n)_n$ is then defined as the enveloping $C^*$-algebra of $(A',p)$. It is important to notice that the direct limit is not unique, in the sense that it highly depends on the choices of maps $\psi_n$. We now state the following proposition that characterizes the inductive limit $A$ in terms of the algebras $A_n$.
\begin{mprop}[\cite{murphy_c-algebras_1990}] Let $(A_n,\psi_n)_n$ be an inductive sequence in the category of $C^*$-algebras. Then there exists an inductive limit $(A,\psi_{n,\infty})$ which satisfies the following:
\begin{itemize}
\item[(i)] $A=\overline{\bigcup_{n\in \mathbb{N}}\psi_{n,\infty}(A_n)}$;
\item[(ii)] For any $n\in \mathbb{N}$ and $a\in A_n$, $\|\psi_{n,\infty}(a_n)\|=\lim_{p\rightarrow\infty}\|\psi_{n,p}(a)\|$.
\item[(ii)] For any $n\in \mathbb{N}$, $a\in \ker\psi_{n,\infty}$ if and only if $\lim_{p\rightarrow\infty}\|\psi_{n,p}(a)\|=0$.
\end{itemize}
\end{mprop}
We consider now the inverse system $\lbrace X_n, \mathbb{N}, \phi_{m,n}\rbrace$ defined in Section \ref{inverse}. To each poset $X_n$, we associate a $C^*$-algebra $(A_n,id_n)$ through the Behncke-Leptin construction. We have then the following identification:
\begin{equation*}
Spec(A_n)\simeq X_n \quad \forall n\in \mathbb{N}.
\end{equation*}
Moreover, using Proposition \ref{pullback2}, the map $\phi_{n,n+1}:X_{n+1}\rightarrow X_n$ induces a pullback map $\phi^*_{n,n+1}:A(X_n)\rightarrow A(X_{n+1})$ for all $n\in \mathbb{N}$. We then have the following diagram in Figure \ref{limitdiagram}.
\begin{figure}[ht!]
\begin{center}
\begin{tikzcd}[cramped, column sep=large, row sep=huge]
A_1 \arrow["id_1",d] \arrow[rightarrow,"\phi^*_{12}",r ] & A_2  \arrow["id_2", d] \arrow[rightarrow,"\phi^*_{23}",r ] & A_3  \arrow["id_3", d]\arrow[rightarrow,r ] & \cdots \arrow[rightarrow,r ] & A_{\infty} \arrow["id", dashed, d]\\
X_1 \arrow[leftarrow,"\phi_{12}",r ] & X_2 \arrow[leftarrow,"\phi_{23}",r ] & X_3 \arrow[leftarrow,r ] & \cdots \arrow[leftarrow,r ] & X_{\infty}
\end{tikzcd}
\end{center}
\caption{Direct system of $C^*$-algebras }
\label{limitdiagram}
\end{figure}
\begin{mprop}The system $\lbrace A_n,\mathbb{N}, \phi^*_{m,n}\rbrace$ forms a direct system.
\end{mprop}
\begin{proof}
We start by recalling that the maps $\phi_{m,n}$ satisfy the coherence properties:
\begin{equation*}
\phi_{l,m} \circ \phi_{m,n} = \phi_{l,n}, \quad l\leq m \leq n, \quad \phi_{n,n}= id_n \quad \forall n\in \mathbb{N}.
\end{equation*}
From this, it follows that for any $l\leq m \leq n$, the following equalities hold:
\begin{align*}
(\Phi_{l,m} \circ \Phi_{m,n} )^{-1} &:= \lbrace y \in (\phi_{l,m} \circ \phi_{m,n})^{-1}(x) : \dim(\sigma'_y)=\dim(\sigma_x) \rbrace,\\
&= \lbrace y \in \phi_{m,n}^{-1} \circ \phi_{l,m}^{-1}(x) : \dim(\sigma'_y)=\dim(\sigma_x) \rbrace,\\
&= \Phi_{m,n}^{-1} \circ \Phi_{l,m}^{-1},
\end{align*}
on one hand; and on the other hand
\begin{align*}
(\Phi_{l,m} \circ \Phi_{m,n} )^{-1}&=\lbrace y \in \phi_{l,n}^{-1}(x) : \dim(\sigma'_y)=\dim(\sigma_x) \rbrace,\\
&=\Phi_{l,n}^{-1}.
\end{align*}
This implies by construction that the pullback maps $\phi_{m,n}^*$ also satisfy the coherence properties: 
\begin{equation}
\phi_{m,n}^* \circ \phi_{l,m}^* = \phi_{l,n}^*, \quad l\leq m \leq n, \quad \phi_{n,n}^*= id_n \quad \forall n\in \mathbb{N}.
\end{equation}
and thus $\lbrace A_n,\mathbb{N}, \phi^*_{m,n}\rbrace$ forms a direct system.
\end{proof}
We can now write the direct limit as
\begin{equation}
A_\infty:=\lim_{\rightarrow} \ (A_n,\phi^*_{n,n+1})_{n\in \mathbb{N}}.
\end{equation}
Let $Z(A_\infty)$ be the center of $A_\infty$; consider the space $MZ_A$ being the space of maximal ideal in $Z(A_\infty)$ equipped with the hull-kernel topology.
From the Gelfand-Naimark theorem \cite[Thm 2.2.4 p.60]{blackadar_operator_2006}, we deduce immediately that $Z(A_\infty)$ is $^{*}$-isomorphic to the space of continuous functions  $C(MZ_A,\mathbb{C})$. Therefore, to prove that $Z(A_\infty)$ is isomorphic to the space of functions $C(M,\mathbb{C})$ over the manifold $M$, we only need to prove that the spaces $\mathfrak{M}$ and $ MZ_A$ are homeomorphic. In fact, we can prove a stronger result:
\begin{mth}\label{main1}
The spectrum $Spec(A_\infty)$ equipped with the hull-kernel topology is homeomorphic to the space $X_\infty$ and 
\begin{equation}
\lim_{\leftarrow}Spec(A_i) \simeq Spec(\lim_{\rightarrow}A_i).
\label{maxid} 
\end{equation} 
\end{mth}
Before proving this result, we recall the definition of a \textit{state} and the interplay with representations. A state $\varphi$ is a positive linear functional with $\varphi(1)=\|\varphi\|=1$. We denote by $S(A)$ the space of states over the $C^*$-algebra $A$ equipped with the weak$^{*}$ topology. In addition, the set $S(A)$ is convex; an \textit{extreme point} of $S(A)$ is called a \textit{pure state} and the set of pure states is denoted by $P(A)$. We will denote the set of extreme points of a convex set $C$ by $\text{ext}(C)$.\\
The GNS construction (see for instance \cite[pp.114-115]{blackadar_operator_2006}) give a one-to-one correspondence between positive linear functional $\varphi$ and (cyclic) representation $(H_\varphi,\pi_\varphi,\xi_\varphi)$.\par
Now let $x\in X_\infty$, then identifying $X_i$ with $Spec(A_i)$, the corresponding representation $\pi_x$ defines a coherent sequence
\begin{equation*}
\pi_x=(\pi_1,\pi_2,\cdots)\in \prod_{i\in \mathbb{N}}Spec(A_i), \quad \text{such that $\pi_m = \phi_{m,n}(\pi_n), $ $\forall m\leq n$}.
\end{equation*}
Moreover, according to the GNS construction, we can associate a pure state $\varphi$ to any irreducible representation $\pi$. Therefore, we have the following coherent sequence of pure states:
\begin{equation}
\varphi_x=(\varphi_1,\varphi_2,\cdots)\in \prod_{i\in \mathbb{N}}S(A_i),
\end{equation}
such that,
\begin{align}
&\varphi_m = \phi_{n,m}(\varphi_n), \label{coh1} \\
&\phi_{l,m}= \phi_{l,n}\circ\phi_{n,m}, \quad \text{if $l\leq m \leq n$ }.
\label{coh2}
\end{align}
Hence, the inverse system of posets $\lbrace X_n,\mathbb{N},\phi_{m,n} \rbrace$ induces an inverse system of states $\lbrace S(A_n),\mathbb{N},\phi_{m,n} \rbrace$.
\begin{mlem}
The inverse limit system $\lbrace S(A_n),\mathbb{N},\phi_{m,n} \rbrace$ is homeomorphic to $S(A_\infty)$.
\label{state}
\end{mlem}
\begin{proof}
For $x\in X_\infty$, the map $\varphi_x$ defines a bounded linear functional on the algebraic inductive limit $A'$ and uniquely extend over $A_\infty$ such that $\| \varphi_x\|=1$. Hence, $\varphi_x\in S(A_\infty)$.\\
Reciprocally, any state $\phi\in S(A_\infty)$ define a state $\varphi_n\in S(A_n)$ defined as follows 
\begin{equation}
\varphi_n := \varphi\circ\phi^*_{n,\infty}(a)
\end{equation}
for any $n\in \mathbb{N}$. In addition, the sequence $(\varphi_n)$ is a coherent sequence satisfying \eqref{coh1} and \eqref{coh2}. Thus, there is a bijection between $\lim_{\leftarrow}S(A_i)$ and $S(A_\infty)$.\par
Finally, the weak$^*$-topology on $S(A_\infty)$ is equivalent to the subspace topology on $\lim_{\leftarrow}S(A_i)$ induced by the product topology on $\prod_{i\in \mathbb{N}}S(A_i)$. This gives us the expected homeomorphism.
\end{proof}
\begin{mlem}
The inverse limit system $\lbrace P(A_n),\mathbb{N},\phi_{m,n} \rbrace$ is homeomorphic to $P(A_\infty)$.
\label{pstate}
\end{mlem}
\begin{proof}
We start by recalling that the inverse limit of convex spaces is convex (this follows from the fact that an arbitrary Cartesian product of convex sets is convex). Therefore, the set $\lim_{\leftarrow}S(A_i)$ is convex. In addition, the set of extreme points of $S(A_i)$ is exactly the set of pure states $P(A_i)$. Using a classical result in convex analysis \cite[Thm.3 p.502]{kadets_course_2018}, the set of extreme points in the product is given by:
\begin{equation}
\text{ext}\left( \prod_{i\in \mathbb{N}}S(A_i)\right) = \prod_{i\in \mathbb{N}}P(A_i)
\end{equation}
Hence, the pure states of $\lim_{\leftarrow}S(A_i)$ is given by the coherent sequences in $\prod_{i\in \mathbb{N}}S(A_i)$ i.e. $\lim_{\leftarrow}P(A_i)$. Similarly, the set of pure states on $A_\infty$ is denoted by $P(A_\infty)$. Consequently, using Lemma \eqref{state} we deduce that
\begin{equation}
\text{ext}\left(\lim_{\leftarrow}S(A_i) \right) = \text{ext}\left( S(A_\infty)\right)= P(A_\infty).
\end{equation}
Finally, we recall that a sequence of states $(\varphi_n)$ on $A_\infty$ converges to a state $\varphi$ in the usual weak topology if and only if the coordinate sequence $(\varphi_n^{i})$ on $A_i$ converges for every $i\in\mathbb{N}$. Therefore, the space $P(A_\infty)$ is homeomorphic to the closed subspace of all systems satisfying \eqref{coh1} in the product space $\prod_{i\in \mathbb{N}}P(A_i)$ i.e
\begin{equation}
\lim_{\leftarrow}P(A_i) \simeq P(A_\infty).
\end{equation}
\end{proof}
\begin{proof}[Proof of Theorem~\ref{maxid}]
Let $\pi\in Spec(A_\infty)$, then by the GNS construction, we can associate to it a pure state $\varphi\in P(A_\infty)$. Using Lemma \ref{pstate}, $\varphi$ in turn correspond to a sequence of pure states $(\varphi_i)$ in $\lim_{\leftarrow}P(A_i)$. Again by Lemma \ref{pstate} and the GNS construction, we associate to $(\varphi_i)$ a coherent sequence in $X_\infty$.\\
Reciprocally, a coherent sequence of irreducible representations in $X_\infty$ correspond to an element in $\lim_{\leftarrow}P(A_i)$ through the GNS construction.\\
Therefore, we can identify $X_\infty$ with $Spec(A_\infty)$ as posets. The homeomorphism follows from the fact that the order topology on $X_\infty$ is equivalent to the hull-kernel topology using the isomorphism \eqref{openideal}.
\end{proof}
\begin{mcor} The sets $\mathfrak{M}$ and $ MZ_A$ are homeomorphic.
\end{mcor}
\begin{proof}
This follows again from the isomorphism \eqref{openideal} where the maximal points in $X_\infty$ correspond to maximal ideals in $Spec(A_\infty)$. Then $\mathfrak{M}$ and $ MZ_A$ are homeomorphic with the subspace topology.
\end{proof}
We have then proven that the $C^*$-algebra $A_\infty$ contains the algebra of continuous functions $C(M,\mathbb{C})$ as its centre. In fact, one can go further in the characterization of the inductive limit using the following result.
\begin{mth}[Dauns-Hofmann {\cite[p.272]{doran_characterizations_1986},\cite{ forger_locally_2013}}] Let $A$ be a unital $C^*$-algebra with centre $Z(A)$. Let $MZ_A$ be the space of maximal ideals of the center $Z(A)$ equipped with the hull-kernel topology. Then $A$ is isometrically $^{*}$-isomorphic to the $C^*$-algebra of all continuous sections $\Gamma (MZ_A,A)$ of the $C^*$-bundle $(\textup{A},\Psi, MZ_A)$ over $MZ_A$. The fibre (stalk) above $x\in MZ_A$ is given by the quotient $\textup{A}_x\simeq \bigslant{A}{xA}$, the isometric $^{*}$-isomorphism is Gelfand's representation $a\mapsto \hat{a}$:
\[
\left\lbrace 
\begin{array}{ccc}
A & \rightarrow & \Gamma (MZ_A,A) \\ 
a & \mapsto & x \mapsto \hat{a}(x)= a +xA \\ 
\end{array} 
\right. 
\]
with $\|\hat{a}\|=\sup_{x\in MZ_A}\|\hat{a}(x)\|$.
\end{mth}
According to the Dauns-Hofmann theorem, the algebra $A_\infty$ is isomorphic to the $C^*$-algebra of continuous sections $\Gamma(M,A_\infty)$ of a $C^*$-bundle $(\textup{A},\Psi, M)$ over the manifold $M$. From the Behncke-Leptin construction, we get the following general form for a section at a point $x\in M$.
\begin{equation}
\hat{a}(x)=\sum_{i\in I_x}\lambda_i(x)\otimes k_i(x) + xA_\infty
\end{equation}
where $I_x$ is a finite indexing set. We see that the central elements are then given by functions $x\mapsto \lambda(x)$ on $M$.\\

We go back now to the commutative subalgebra $\mathfrak{A}$ defined in Equation \eqref{commutative} and show how it can be used to approximate $C(M)$. In the rest of this work, we will identify $C(M)$ with the centre $Z(A_\infty)$ and denote by $\mathfrak{A}_n$ the commutative subalgebra in $A(X_n)$.\\
\begin{mprop} The space of continuous function $C(M)$ is approximated by the system of commutative subalgebras $(\mathfrak{A}_n,\phi^*_{n,\infty})$ in the following sense:
\begin{equation}
C(M)=\overline{\bigcup_{n\in \mathbb{N}}\phi^*_{n,\infty}(\mathfrak{A}_n)} \cap C(M).
\end{equation}
\label{main2}
\end{mprop}
\begin{proof}
First, let us recall that, by Axiom 5), an element $a_n\in \mathfrak{A}_n$ is determined a map $\hat{a}_n:X_n\rightarrow A_n$ such 
\begin{equation}
\hat{a}_n(x)=\sum_{i\in I_x}\lambda_i(x)1_{H(x)}.
\end{equation}
When restricted to the set of maximal points $\mathfrak{M}_n$, $a_n$ acts as a scalar:
\begin{equation}
\pi_x(a_n)=\lambda_x
\end{equation}
where $\lambda_x\in \mathbb{C}$. Then, using the map $\phi_{n,\infty}$, we notice that $a_n$ defines a piecewise-linear function on $M$:
\begin{equation}
\phi^*_{n,\infty}(a_n)= a_n\circ \phi_{n,\infty}:M\rightarrow\mathbb{C}, \quad
a_n\circ \phi_{n,\infty}(y)=\lambda_{\phi_{n,\infty}(y)}.
\end{equation}
Therefore, any continuous function $g\in C(M)$ can be uniformly approximated arbitrarily closely by a function of the form $a_n\circ \phi_{n,\infty}$, for some sufficiently large $n$.
\end{proof}
Finally, using the smooth structure, we can define the subalgebras 
\begin{equation}
Z^{k}(A_\infty):=C^{k}(M)
\end{equation}
of $k$-differentiable functions. In the rest of this work, we will focus on the subalgebra $Z^{\infty}(A_\infty)$ and its approximation given by the equality:
\begin{equation}
Z^{\infty}(A_\infty)=\overline{\bigcup_{n\in \mathbb{N}}\phi^*_{n,\infty}(\mathfrak{A}_n) \cap Z^{\infty}(A_\infty)}.
\end{equation}
\subsubsection{Direct limit of representations}
Similarly, we associate a representation space $H(X_n)$ (defined in Equation \eqref{representation}) to every space $X_n$. Moreover, a continuous surjection $\phi:X'\rightarrow X$ between posets induces a isometry $\psi:H(X)\rightarrow H(X')$ between representations. The construction of $\psi$ follows mutatis mutandis the same steps that the one of $\phi^*$; therefore, we will keep the same notations and directly states the results. We define $\psi:H(X)\rightarrow H(X')$ as follows :
\begin{equation}
H(X)=\oplus_{x\in X}H(x), \quad \psi\left( \oplus_{x\in X} \xi_x\right) = \oplus_{y\in X'} \xi_y
\end{equation}
where,
\begin{equation}
\xi_y = \left\lbrace 
\begin{array}{cc}
\xi_\phi(y) & \text{if $\dim(\sigma_y)=\dim(\sigma_{\phi(y)})$,} \\ 
0 & \text{otherwise}.
\end{array} 
\right. 
\label{hilbert}
\end{equation}
Therefore, the inverse system of posets $\lbrace X_n,\mathbb{N},\phi_{m,n} \rbrace$ induces a direct system of Hilbert spaces $\lbrace H_n,\mathbb{N},\psi_{m,n} \rbrace$, where $H_n$ denotes the Hilbert space $H(X_n)$.
\begin{mprop}The system $\lbrace H_n,\mathbb{N}, \psi_{m,n}\rbrace$ forms a direct system.
\end{mprop}
Hence, we can construct the direct limit of representations $(H_n,\psi_n)$ as a subspace of the direct sum:
\begin{equation}
\bigoplus_{n\in \mathbb{N}}H_n=\left\lbrace (h_n)_{n\in \mathbb{N}} : h_n\in H_n, \sum_{n=1}^\infty \|h_n\|^2_{H_n}<\infty\right\rbrace
\end{equation}
equipped with an inner product $\left\langle .,.\right\rangle $ given by:
\begin{equation}
\left\langle g,h\right\rangle=\sum_{n=1}^\infty \left\langle g_n,h_n\right\rangle_{H_n}.
\end{equation}
The algebraic direct limit is defined as 
\begin{equation}
H'=\left\lbrace (h_n)_{n\in \mathbb{N}} \in \bigoplus_{n\in \mathbb{N}}H_n : \psi_{n,n+1}(h_n)=h_{n+1} \right\rbrace .
\end{equation}
The resulting Hilbert space is obtained from the closure of $H'$ and will be denoted by 
\begin{equation}
H_{\infty}:=\lim_{\rightarrow}(H_n,\psi_n)_{n\in \mathbb{N}}.
\end{equation}
The direct system $\lbrace H_n,\mathbb{N}, \psi_{m,n}\rbrace$ induces a direct system on the irreducible representations $\lbrace \mathcal{H}_n,\mathbb{N}, \psi_{m,n}\rbrace$; we denote the limit $\mathcal{H}_\infty$. We have then the following characterization of this limit space.
\begin{mth} The Hilbert space $L^2(M)$ of square integrable functions over the manifold $M$ is a subspace of $\mathcal{H}_\infty$:
\begin{equation}
\mathcal{H}_\infty = L^2(M)\oplus \mathcal{H}_\omega.
\end{equation}
\end{mth}
\begin{proof} Using Axiom 5 in the Behncke-Leptin construction, for any $n\in \mathbb{N}$, we have the following decomposition:
\begin{equation}
\mathcal{H}_n=\bigoplus_{x\in \mathfrak{M}_n}\mathcal{H}^n_x\oplus \bigoplus_{x\in \mathfrak{M}_n^c}\mathcal{H}^n_x\simeq \mathbb{C}^{|\mathfrak{M}_n|} \oplus \bigoplus_{x\in \mathfrak{M}_n^c}\mathcal{H}^n_x.
\end{equation}
Now, let us recall that the commutative subalgebra $\mathfrak{A}_n$ given by
\begin{equation}
\mathfrak{A}_n=\oplus_{x\in \mathfrak{M}_n}1_{H(x)}, \quad a=\sum_{x\in \mathfrak{M}_n}\lambda(x)1_{H(x)}
\end{equation}
is completely determined by the representation $(\mathbb{C}^{|\mathfrak{M}_n|},\oplus_{x\in \mathfrak{M}_n} \pi_x)$:
\begin{equation}
\oplus_{x\in \mathfrak{M}_n} \pi_x:\mathfrak{A}_n \rightarrow \mathbb{C}^{|\mathfrak{M}_n|}, \quad a \mapsto (\lambda(x_1), \lambda(x_2),\cdots,\lambda(x_{|\mathfrak{M}_n|})).
\end{equation}
Through this isomorphism of vector space, we can identify $\bigoplus_{x\in \mathfrak{M}_n}\mathcal{H}^n_x$ with the image of $\mathfrak{A}_n$ and denote it by $\hat{\mathfrak{A}}_n$:
\begin{equation}
\mathcal{H}_n= \hat{\mathfrak{A}}_n \oplus \bigoplus_{x\in \mathfrak{M}_n^c}\mathcal{H}^n_x.
\end{equation}
Moreover, because $\phi_{n,n+1}(\mathfrak{M}_{n+1})=\mathfrak{M}_n$ then by definition of $\psi_{n,n+1}$, we have:
\begin{align}
\psi_{n,n+1}(\mathcal{H}_n)&=\phi^*_{n,n+1}(\hat{\mathfrak{A}}_{n})\oplus \bigoplus_{x\in \mathfrak{M}_n^c}\psi_{n,n+1}(\mathcal{H}^n_x)
\end{align}
for every $n\in \mathbb{N}$. Therefore, we have for every $n\in \mathbb{N}$:
\begin{equation}
\psi_{n,\infty}(\mathcal{H}_n)=\phi^*_{n,\infty}(\hat{\mathfrak{A}}_{n})\oplus \bigoplus_{x\in \mathfrak{M}_n^c}\psi_{n,n+1}(\mathcal{H}^n_x).
\end{equation}
Hence the direct limit $H_\infty$ decomposes as follows:
\begin{equation}
H_\infty=H\oplus H_\omega,
\end{equation}
where $H_\omega$ is an infinite dimensional Hilbert space and
\begin{equation}
H=\overline{\oplus_n \phi^*_{n,\infty}(\hat{\mathfrak{A}}_{n})}= \overline{ \left\lbrace a\in C(M), \|a\|_{H_\infty}<\infty  \right\rbrace } \equiv L^2(M).
\end{equation}
\end{proof}
\subsubsection{Cubulation: example of a lattice}
We now conclude this section with the specific case of  a $C^*$-algebra over a lattice $\Lambda$ seen as a cubulation of $\mathbb{R}^d$. The lattice $\Lambda$ can be written as a direct product of a line lattice $L$. Hence, we can the algebra $A(\Lambda)$ relate them to the tensor product of algebras $A(L)$ over $L$. First, we need to recall the following result on the structure space of tensor product of $C^*$-algebras.
\begin{mprop}[Wulfsohn {\cite{wulfsohn_primitive_1968}}] Let $A$ and $B$ be separable $C^*$-algebras and $A\otimes B$ their $C^*$-tensor product. The mapping
\label{primtensor1}
\begin{equation*}
\alpha: Prim(A)\times Prim(B)\rightarrow Prim(A\otimes B), \quad \alpha(\mathfrak{a},\mathfrak{b})=\mathfrak{a}\otimes B+A\otimes \mathfrak{b}
\end{equation*}
is a homeomorphism.
\end{mprop}
This result immediately gives us that tensor $C^*$-algebras can be seen as $C^*$-algebras over Cartesian product of posets.
\begin{mcor}Let $X$ and $Y$ be topological spaces. If $(A,\psi_A)$ and $(B,\psi_B)$ are separable $C^*$-algebra over $X$ (respectively over $Y$), then the pair $(A\otimes B,\psi_A\times \psi_B)$ is a separable $C^*$-algebra over $X\times Y$ with the product topology.
\label{primtensor2}
\end{mcor}
Let $\Lambda$ be the $d$-dimensional, we can write it as the direct product of $d$ line lattices:
\begin{equation*}
\Lambda = L\times \cdots \times L.
\end{equation*}
Let $(A(L),\psi_L)$ be a $C^*$-algebra over $L$. Then using \ref{primtensor1} and \ref{primtensor2} we can associate the $C^*$-algebra over $\Lambda$:
\begin{equation}
A(\Lambda)=A(L)\otimes \cdots \otimes A(L), \quad \psi_\Lambda = \mathrm{\Pi} \psi_L.
\end{equation}
 Similarly to the previous section, we construct a sequence of refined lattice $(\Lambda_n,\pi_{n})$ and construct the direct limits of $C^*$-algebras $(A(\Lambda_n),\pi^*_{n})$ with their representations $(H_n,\psi_{n})$. We can then directly state the following result, which a special case when $M=\mathbb{R}^d$.
\begin{mcor}
The centre of the limit $C^*$-algebra $A_\infty$, $Z(A_\infty)$ is isometrically $^{*}$-isomorphic to $\mathcal{C}(\mathbb{R}^n)$ acting on $L^2(\mathbb{R}^n)$ as a subspace of $H_\infty$.
\end{mcor}
\section{Geometry over a triangulation}
\label{Gem}
The last piece remaining to define, in order to complete this triptych, is the differential geometry. This will be done using the machinery of noncommutative differential geometry, as explained in the introduction.
\subsection{Finite spectral triple}
Let $(A(X),id)$ be a $C^*$-algebra over a poset $X$ induced by a triangulation of a compact Riemannian manifold $(M,g)$ of dimension $d$. \\
We will denote by $\mathfrak{M}$ the set of maximal points in $X$ and by $\mathfrak{A}$ the commutative subalgebra of $A$ defined by Equation \eqref{commutative}. We then immediately notice that
\begin{equation}
\mathfrak{h}=\bigoplus_{x\in \mathfrak{M}}\mathbb{C},\quad \pi=\bigoplus_{x\in \mathfrak{M}}\pi_x.
\end{equation}
defines a faithful representation of $\mathfrak{A}$.\\
Consider now the pair $(\mathfrak{h},\mathfrak{h}^*)$ where $\mathfrak{h}$ and $\mathfrak{h}^*$ have both dimension $m$. Define the even dimensional representation of $\mathfrak{A}$
\begin{equation}
\mathfrak{H}(X):= \mathfrak{h}\oplus \mathfrak{h}^*, \quad \rho=\pi\oplus \pi^*
\label{splitting}
\end{equation}
where the adjoint representation is given by $\pi^*(a)=-\pi^t(a)$ for any $a\in \mathfrak{A}$. The triple $(\mathfrak{A},\mathfrak{H},\rho)$ embeds the commutative algebra $\mathfrak{A}$ into the \textit{Cartan subalgebra} $\mathsf{h}$ of the Lie algebra $\mathfrak{gl}(2m,\mathbb{C})$. \\
The space of bounded operators $B(\mathfrak{H})$ can be identified with $M_{2m}(\mathbb{C})$. We define the parity element $\gamma\in M_{2m}(\mathbb{C})$ such that 
\begin{equation}
\gamma =
\left( 
\begin{array}{cc}
1_m & 0 \\
0 & -1_m
\end{array} 
\right) 
\end{equation}
where the eigenspace decomposition correspond to the splitting \eqref{splitting}. This in turns defines a $\mathbb{Z}_2$-grading on $M_{2m}(\mathbb{C})$. The space $M_{2m}$ can be accordingly written as a direct sum
\begin{equation}
M_{2m}=M_{2m}^+\oplus M_{2m}^-
\end{equation}
of even and odd elements, where $a\in M_{2m}$ is even if it commutes with $\gamma$ and odd if it anticommutes. In fact, even elements will correspond to block diagonal elements and odd elements to off-diagonal with respect to the representation space $\mathfrak{H}$. Under this grading, the algebra $\mathfrak{A}$ is represented as the subspace of diagonal matrices, i.e.
\begin{equation}
\mathfrak{A}\overset{\rho}{\longrightarrow}\mathsf{h}\longhookrightarrow M_{2m}^+(\mathbb{C}).
\end{equation}
\begin{mrmk}
The data $(\mathfrak{A},\mathfrak{H},\pi)$ can also be localized to an open set $U\subset X$. Consider the restriction functor $r_X^U$ (see \cite{meyer_c-algebras_2008}) and define the restriction $A(U):=(r_X^UA)$ of $A$ to the open set $U$. 
Similarly, $\mathfrak{A}(U)$ defines a restriction of $\mathfrak{A}$ to $U$. Let $\mathfrak{M}_U$ the subset of $\mathfrak{M}$ of maximal points in $U$. Again, we have that
\begin{equation*}
\mathfrak{H}_U=\bigoplus_{x\in \mathfrak{M}_U}\mathbb{C},\quad \pi|_U=\bigoplus_{x\in \mathfrak{M}_U}\pi_x,
\end{equation*}
is a representation of $\mathfrak{A}(U)$.
\end{mrmk}
\begin{mdef}[Spectral triple] A \textit{spectral triple} is the data $(\mathcal{A},\mathcal{H},D)$ where:
\begin{itemize}
\item[(i)] $\mathcal{A}$ is a real or complex $*$-algebra;
\item[(ii)] $\mathcal{H}$ is a Hilbert space and a left-representation $(\pi,\mathcal{H})$ of $A$ in $\mathcal{B}(\mathcal{H})$;
\item[(iii)] $D$ is a \textit{Dirac operator}, which is a self-adjoint operator on $\mathcal{H}$.
\end{itemize}
If in addition, $\mathcal{H}$ is equipped with a $\mathbb{Z}_2$-grading i.e. there exists a unitary self-adjoint operator $\gamma\in \mathcal{B}(\mathcal{H})$ such that
\begin{itemize}
\item[1)] $\left[\gamma,\pi(a)\right]=0$ \ \text{for all $a\in \mathcal{A}$}, 
\item[2)] $\gamma$ anticommutes with $D$,
\end{itemize}
then the spectral triple is said to be \textit{even}. Otherwise, it is said to be \textit{odd}. In the case where $\mathcal{H}$ is finite dimensional, then the triple $(\mathcal{A},\mathcal{H},D)$ is called a \text{discrete spectral triple}.
\end{mdef}
We consider the finite dimensional algebra $(\mathfrak{A},\mathfrak{H})$ a Dirac operator $D$ chosen as an odd element of $M_{2m}(\mathbb{C})$ of the form
\begin{equation}
D=\frac{i}{h}
\left( 
\begin{array}{cc}
0 &  D^-\\
D^+ & 0
\end{array} 
\right) 
\end{equation}
where $D^+, D^-\in M_{2m}(\mathbb{R})$ and satisfy $D^-=-(D^+)^*$. We then form the finite spectral triple $(\mathfrak{A},\mathfrak{H},D)$; this triple is even with the grading induced by $\gamma$. \par
Using this structure, we can then define a graded derivation $da$ for $a\in M_{2m}(\mathbb{C})$ through a graded commutator,
\begin{equation}
da=-\left[D,a\right]:=Da - \epsilon_aaD
\label{diff}
\end{equation}
where $\epsilon_a=1$ if a is even and $\epsilon_a=-1$ if a is odd. Using the representation $\rho$, it also induces a derivation on $\mathfrak{A}$. Furthermore, notice that the derivative $d$ coincides (modulo the grading) with the adjoint operator $ad_D$. We can then study the differential structure on $\mathfrak{A}$ by identifying $M_{2m}(\mathbb{C})$ as the Lie algebra $\mathfrak{gl}_{2m}(\mathbb{C})$ with Cartan subalgebra $\mathsf{h}$. For convenience, we then equip $\mathsf{h}$ with the inner product:
\begin{equation}
\left\langle h, h' \right\rangle := Tr(h^*h').
\end{equation} 
We can then identify $\mathsf{h}$ with its dual $\mathsf{h}^*$ i.e. the set of linear functionals acting on $\mathsf{h}$. Now, recall that a nonzero element $\alpha\in \mathsf{h}$ is a \textit{root} of $\mathfrak{gl}_n(\mathbb{C})$ relative to $\mathsf{h}$ if there exists a nonzero $x\in \mathfrak{gl}_n(\mathbb{C})$ such that
\begin{equation}
[x,h]=\alpha(h)x,
\end{equation}
for all $h\in \mathsf{h}$. In particular, the standard matrix basis elements $e_{ij}$ satisfies $he_{ij}=\lambda_ie_{ij}$ and $e_{ij}h=\lambda_je_{ij}$ for all $h\in \mathsf{h}$. Thus, 
 \begin{equation}
 [h,e_{ij}]=(\lambda_j-\lambda_i)h,
 \end{equation}
 showing that $e_{ij}$ are simultaneous eigenvectors for $ad_h$. Now let $a\in \mathfrak{A}$ be described as an element of $\mathsf{h}$ through the representation $\rho$:
\begin{equation}
\rho(a)=\left( 
\begin{array}{cccccc}
\lambda_1 &  &  &  &  &  \\ 
 & \ddots &  &  &  &  \\ 
 &  & \lambda_m &  &  &  \\ 
 &  &  & \lambda_1 &  &  \\ 
 &  &  &  & \ddots &  \\ 
 &  &  &  &  &\lambda_m \
\end{array} 
\right).
\label{cartan}
\end{equation}
Following the definition, we can write the operator $D$ as a linear combination of elements $e_{ij}$:
\begin{equation}
D=\sum_{i<j}\omega_{ij}\hat{e}_{ij}
\end{equation}
where $\hat{e}_{ij}=e_{ij} - e_{ji}$. Then the derivation $d$ acts on an element $a\in \mathfrak{A}$ as
\begin{equation}
da=\sum_{i<j} \omega_{ij}\alpha_{ij}(a)\hat{e}_{ij}
\label{da}
\end{equation}
where the roots are given by $\alpha_{ij}(a)=\lambda_j - \lambda_i$.
\subsubsection{Graded differential algebra}
It is possible to construct over $M_{2m}^+$ a $\mathbb{N}$-graded differential algebra $\Omega_D^*=\Omega_D^*(M_{2m}^+)$ based on formula \eqref{diff}. Define $\Omega_D^0=M_{2m}^+$ and let 
\begin{equation}
\Omega_D^1=d\Omega_D^0 \subset M_{2m}^-
\end{equation}
be the $M_{2m}^+$-module generated by the image of $\Omega_D^0$ in $M_{2m}^-$ under $d$. Then for each $p$, we let $\Ima d^2$ be the submodule of $d\Omega_D^{p-1}$ consisting of those elements which contain a factor which is the image of $d^2$ and define 
\begin{equation}
\Omega^p_D=d\Omega_D^{p-1}/\Ima d^2.
\end{equation}
Therefore since $\Omega_D^p\cdot \Omega_D^q\subset \Omega_D^{p+q}$ the complex $\Omega_D^*$ define as 
\begin{equation}
\Omega_D^*=\bigoplus_{p\geq 0}\Omega_D^p
\end{equation}
is a differential graded algebra. The $\Omega_D^p$ need not vanish for large values of $p$. In addition, it follows by construction that the $\Omega_D^p$ are generated by the $da$ as follows
\begin{equation}
\Omega_D^p = \left\lbrace a_0da_1\cdots da_n, \ a_i\in M_{2m}^+(\mathbb{C}) \ \forall i\right\rbrace.
\end{equation}
However, we would like to restrict to elements in $\mathfrak{A}$ seen as a subset of $M^{+}_{2m}$ through the representation $\rho$. We then define $\Omega^*_D(\mathfrak{A})$ in the exact same. In particular, we have 
\begin{equation}
\Omega_D^1(\mathfrak{A}) = \left\lbrace a_0da_1, \ a_i\in \mathfrak{A}, \ i=1,2\right\rbrace.
\end{equation}
We define an inner product on $B(\mathfrak{H})$ given by
\begin{equation}
(A,B)_{B(\mathfrak{H})}=Tr(B^*A),
\end{equation}
and inducing a Hilbert space structure on $\Omega^k_D(\mathfrak{A})$ for any $k$.
\subsubsection{Laplace operator}
Following the definitions, we see that the differential $da$ is not an element of $\mathfrak{A}$ in general, but is in $B(\mathfrak{H})$ nonetheless. Let $p$ be the orthogonal projection operator on $\rho(\mathfrak{A})$ with respect to this inner product:
\begin{equation}
B(\mathfrak{H}) = \rho(\mathfrak{A})\oplus \rho(\mathfrak{A})^{\perp}.
\end{equation}
We can now introduce the adjoint operator $\delta : B(\mathfrak{H})\rightarrow \rho(\mathfrak{A})$ of the differential $d$ using Riesz representation theorem
\begin{equation}
(b,da)_{B(\mathfrak{H})}=(a',a)_{\rho(\mathfrak{A})}
\end{equation} 
and we set $\delta b := a'$.
\begin{mprop} The adjoint map $\delta$  to the derivation $d:\mathfrak{A}\rightarrow B(\mathfrak{H})$ is given by
\begin{equation*}
\delta:B(\mathfrak{H})\rightarrow \mathfrak{A}, \quad \delta(b) = p[D,b].
\end{equation*}
\end{mprop}
\begin{proof}
Using the fact that $D$ is hermitian, we first have
\begin{equation}
(b,da)_{B(\mathfrak{H})} = (b,[D,a])_{B(\mathfrak{H})} =  ([D,b],a)_{B(\mathfrak{H})}.
\end{equation}
Then, since $a\in \rho(\mathfrak{A})$ and $p^*=p$, it follows that:
\begin{equation}
([D,b],a)_{B(\mathfrak{H})} = ([D,b],pa)_{B(\mathfrak{H})} = (p[D,b],a)_{\rho(\mathfrak{A})}. 
\end{equation}
\end{proof}
It is then straightforward to define a Laplace operator on $\mathfrak{A}$.
\begin{mdef}(Laplacian) The Laplace operator $\Delta$ is given by:
\begin{equation*}
\Delta: \mathfrak{A}\rightarrow \mathfrak{A},\quad \Delta(a) := -\delta d  a = - p\left[ D,\left[ D,b\right]\right],
\end{equation*}
where $p$ is the orthogonal projection on $\mathfrak{A}$.
\end{mdef}
We can now state and prove a Hodge-like decomposition on $\Omega^*_D(\mathfrak{A})$.
\begin{mprop}[Hodge-de Rham decomposition] The Laplacian $\Delta$ on $\Omega_D(\mathfrak{A})$ satisfies the following properties: 
\begin{itemize}
\item[i)] $\Delta\geq 0$ in the Hilbert space $(\Omega_D(\mathfrak{A}),(\cdot,\cdot))$, 
\item[ii)] $\Delta\alpha =0$ if and only if $d\alpha=0$,
\item[iii)]$\mathfrak{A}=\delta \Omega_D(\mathfrak{A})\oplus \ker(\Delta)$ is an orthogonal decomposition of $\Omega_D(\mathfrak{A})$ with respect to $(\cdot,\cdot)$. 
\end{itemize}
In addition, we will call harmonic these elements $\alpha\in \ker(\Delta)$.
\end{mprop}
\begin{proof}
Let $a\in \mathfrak{A}$, then, $( \Delta(a),a) =\|da\|^2$ which proves $i)$. The inclusion $\ker(\Delta)\subseteq \ker(d)$ follows from $i)$; the inclusion $\ker(d)\subseteq \ker(\Delta)$ is immediate from the definition of $\Delta$. This proves $ii)$. Finally, since $\delta$ is the adjoint to $d$, we have the following decomposition in finite dimension
\begin{equation}
\mathfrak{A}=\ker(d)\oplus \delta \Omega_D(\mathfrak{A});
\end{equation}
thus, $iii)$ follows from $ii)$.
\end{proof} 
\subsection{Dirac operator associated to a graph}
So far, we have worked with a generic Dirac operator $D$, the only restriction being that $D$ has to be hermitian and odd according to the grading. Nevertheless, one can exhibit a deeper connection between the space $X$ (or equivalently the spectrum of $\mathfrak{A}$) and the Dirac operator. We first need to restrict the space of admissible Dirac operators.
\begin{mdef}[Admissible Dirac operators]
Let $D\in M_{2m}(\mathbb{C})$ be an odd and hermitian matrix and let $\omega_{ij}$ be the coefficients of the block $D^-$. We say that $D$ is an admissible Dirac operator associate to $X$ if it satisfies the additional condition:
\begin{align*}
a) \ &\text{vertices $i$ and $j$ do not share an edge} \Leftrightarrow \omega_{ij}=0, \ \forall i,j\in \mathfrak{M},\\
b) \ & \text{the eigenvalues $\mu_n$ satisfy the asymptotic }\ \mu_n(D)=O( h^{-1}).
\end{align*}
We denote by $\mathcal{D}(X)$ the set of all admissible Dirac operators and by $\mathcal{D}_{\mathbb{R}}(X)$ the set of real admissible Dirac operators.
\end{mdef}
\begin{mex} The prototypical example is given by the \textit{combinatorial Dirac operator,} for which:
\begin{equation*}
\omega_{ij}:=\left\lbrace 
\begin{array}{cc}
1 & \text{if $i \sim j$}, \\ 
0 & \text{otherwise.}
\end{array} 
\right. 
\end{equation*}
\end{mex}
\begin{mprop} The graded algebra $\Omega_D(M_{2m}^+(\mathbb{C}))$ is invariant under the change $D\mapsto D'$ in $\mathcal{D}(X)$ i.e. 
\begin{equation*}
\Omega_D(M_{2m}^+(\mathbb{C}))=\Omega_{D'}(M_{2m}^+(\mathbb{C}))
\end{equation*}
for any $D,D'\in \mathcal{D}(X)$.
\end{mprop}
\begin{proof}
 The algebra $\Omega^p(\mathfrak{A})$ is generated by elements of the form $a_0da_1\cdots da_n$, with  $ a_i\in \mathfrak{A}$ for all $0\leq i \leq n$. Now, we recall that for an element $a\in \mathfrak{A}$, 
\begin{equation}
da=\sum_{i<j} \alpha_{ij}(a)\hat{e}_{ij}.
\end{equation}
Therefore, $\Omega^p(\mathfrak{A})$ is generated by basis elements $\lbrace\hat{e}_{ij}\rbrace$ where an element $\hat{e}_{ij}$ is a generator if and only if vertices $i$ a $j$ share an edge.
\end{proof}
\subsection{A first example: the lattice}
We now come back to the example of a $C^*$-algebra over a line lattice $L$ denoted by $\mathfrak{A}(L)$ before moving to the case of a $d$-dimensional lattice $\Lambda$. For the line lattice, we let the Dirac operator $D$ to be an odd element of $M_{2m}(\mathbb{C})$ of the form:
\begin{equation}
D=\frac{i}{h}
\left( 
\begin{array}{cc}
0 &  D^-\\
D^+ & 0
\end{array} 
\right) 
\end{equation}
with $(D^+)^*=-D^-$ and where $D^-$ is given by
\begin{equation}
D^-=\left( 
\begin{array}{ccccc}
0 & 1 & 0 & \cdots & 0 \\ 
\vdots & \ddots & \ddots & \ddots & \vdots \\ 
\vdots &  & \ddots & \ddots & 0 \\ 
\vdots &  & &\ddots  & 1\\
0 & \cdots & \cdots & \cdots & 0
\end{array} 
\right) .
\end{equation}
\begin{mprop} For every element $da\in \Omega^1_{D}(\mathfrak{A})$, the spectrum $\sigma (da)$ of the operator $da$ is given by:
\begin{equation*}
\sigma(da)=\left\lbrace \pm \frac{1}{h}(\lambda_{j+1} - \lambda_j) : 1\leq j \leq m-1\right\rbrace \cup \lbrace 0\rbrace
\end{equation*}
Moreover, we have the commutativity relation 
\begin{equation}
[da,db]=0, \quad \forall a,b \in \mathfrak{A}.
\end{equation}
\label{Spectrum}
\label{main5}
\end{mprop}
\begin{proof}
Using Equation \eqref{da}, we can write the commutator $da$ as 
\begin{equation}
da=\sum_{j=1}^{m-1}\alpha_{jj+m+1}(a)\hat{e}_{jj+m+1}
\end{equation}
with the roots:
\begin{equation}
\alpha_{jj+m+1}(a) = \frac{1}{h}(\lambda_{j+1} - \lambda_{j}), \quad j\in \lbrace 1, \cdots, m-1\rbrace.
\end{equation}
Then, we notice that the operator $da^*da$ is a diagonal operator with diagonal entries given by:
\begin{equation}
\beta_{jj}=\frac{1}{h^2}(\lambda_{j+1}-\lambda_{j})^2=\beta_{j+m+1j+m+1}, \
\end{equation}
for $j\in \lbrace 1, \cdots, m - 1 \rbrace$ and $\beta_{m,m}=\beta_{m+1,m+1}=0$. Thus, the eigenvalues of $da$ are obtained as the square roots of the previous diagonal coefficients. Finally, the commutativity follows again from the fact that $dadb$ is a diagonal operator. 
\end{proof}
Following the result on the spectrum of $da$, we can deduce a result on the states. In the finite dimensional case where $A=M_{2m}(\mathbb{C})$, a density matrix $\omega$ i.e. an operator with a graded trace $Tr_s(\omega)=1$ defines a state. over $A$. Then, we can introduce the \textit{expectation map} $a\mapsto <a>_\omega$ with respect to $\omega$ such that
\begin{equation}
<a>_\omega = Tr_s(\omega a).
\end{equation}
\begin{mprop} There exists a density matrix $\omega$ with eigenvalues $\lbrace\mu_k\rbrace_k^{2m}$, such that the expectation value is given by
\begin{equation*}
<da>_\omega = Tr_s(\omega da) = \frac{i}{h}\sum_{k=1}^{2m}\mu_k(\lambda_{k+1}-\lambda_k).
\end{equation*}
for any element $da\in \Omega_D^1(\mathfrak{A})$.
\end{mprop}
\begin{proof}
By Proposition \ref{Spectrum}, we know that the algebra generated by $da$ for $a$ running in $\mathfrak{A}$ is commutative. Then, it admits a common spectral decomposition. Therefore, one can choose $\omega=db$, for some element $b\in \mathfrak{A}$ such that the graded trace $Tr_s(db)=1$, which conclude the proof.
\end{proof}
\begin{mrmk} The last proposition is of importance in the approximation of differential operators. Indeed, it is well known that any finite difference formula for the first derivative can be written as a convex combination of the two-points approximation. It follows that the perspective is shifted in this context; instead of looking at the pointwise discretization of a derivative, one can study the density matrix $\omega$. It is our hope that this change of paradigm, together with the machinery of $C^*$-algebra, allows us to produce new results in discretization of differential operators.
\end{mrmk}
\subsubsection{Direct limit of spectral triples}
We can now complete the construction in the case of the lattice. Recall that we have defined a direct system of $C^*$-algebras $(A_n,\phi_{m,n}^*)$ over an inverse system $(L_n,\phi_{m,n})$ of lattices. We can now associate a Dirac operator $D_n$ to each algebra $A_n$. We will work on the infinite collection $\left\lbrace \mathfrak{A}_n : n\in \mathbb{N} \right\rbrace $ of commutative subalgebras $C^*$-algebras. In this case, we have identified each of the $\mathfrak{A}_n$ with the Cartan subalgebras $\mathrm{h}_i$ inside the finite dimensional algebras $B_n=M_{2m_n}(\mathbb{C})$ where $m_n\rightarrow \infty$ when $n\rightarrow \infty$. We can then construct the product 
\begin{equation}
B_\omega= \prod_{n\in \mathbb{N}}B_n=\lbrace (a_n) : \|a_n\|=\sup \|a_n\|<\infty\rbrace .
\end{equation}
Let $a$ be an element in $C^\infty(\mathbb{R})$, then there exists a sequence 
$(a_i)$ such that
\begin{equation}
a=(a_0, a_1, \cdots, a_n, \cdots)\in  \prod_{n\in \mathbb{N}} \mathfrak{A}_n.
\end{equation}
We define a spectral triple on $B_\omega$ by introducing the Dirac operator $D$ as the sequence 
\begin{equation}
D=(D_0,D_1,\cdots, D_n, \cdots ) \in \prod_{n\in \mathbb{N}} M_{2m_n}^-(\mathbb{C}).
\end{equation}
This in turns induces a spectral triple on $ \prod_{n\in \mathbb{N}} \mathfrak{A}_n$ along with the commutator:
\begin{equation}
[D,a]=([D_0,a_0],[D_1,a_1],\cdots, [D_n,a_n], \cdots )\in \prod_{n\in \mathbb{N}} M_{2m_n}^-(\mathbb{C}).
\label{commutator}
\end{equation}
We can then characterize the operator $[D,a]$ and relate it to the classical differential on $C^\infty(M)$.
\begin{mlem}
The spectrum $\sigma_{B_\omega}(x)$ of an element $x=(\cdots, x_n, \cdots)\in B_\omega$ is given by
\begin{equation*}
\sigma_{B_\omega}(x)=\overline{\cup_n\sigma_{B_n}(x_n)}.
\end{equation*}
\end{mlem}
\begin{proof}
$ b =(\cdots, b_n, \cdots )\in B_\omega $ is invertible if and only if each $ b_n$ is invertible and $\lbrace \|b_n\|^{-1} : n\in \mathbb{N}\rbrace $ is bounded. Thus, $\sigma_{B_n}(x_n)\subset \sigma_{B_\omega}(x)$ for all $n\in \mathbb{N}$, Therefore, we have the first inclusion $S:=\overline{\cup_n\sigma_{B_n}(x_n)}\subseteq \sigma_{B_\omega}(x)$. Reciprocally, if $\lambda \in \mathbb{C}\backslash S$, then $x_n-\lambda 1$ is invertible in $B_n$ for each $n$ and $\|(x_n-\lambda 1)^{-1}\|\leq d(\lambda, S)$, where $d(\lambda,S)$ is the distance from $\lambda$ to $S$, therefore 
\begin{equation}
(x-\lambda 1)^{-1}=(\cdots,(x_n-\lambda 1)^{-1},\cdots )\in B_\omega .
\end{equation} 
\end{proof}
Then, for an element $a\in \prod_{n\in \mathbb{N}} \mathfrak{A}_n$, the spectrum $\sigma_{B_\omega}([D,a])$ of the operator $[D,a]$ is given by:
\begin{equation}
\sigma_{B_\omega}([D,a])=\overline{\cup_n\sigma_{B_n}([D_n,a_n])}.
\end{equation}
We restrict now to an element $a\in\prod_{n\in \mathbb{N}} \mathfrak{A}_n \cap C^\infty(\mathbb{R})$; using Proposition \ref{Spectrum}, we have:
\begin{equation}
\sigma_{B_n}([D_n,a_n])=\left( \cdots, \frac{a_{x^n_{j+1}}-a_{x^n_j}}{h_n} ,\cdots\right) =\left( \cdots, \ell(a)(x^n_j) ,\cdots\right) 
\end{equation}
where, $a_{x^n_j}=a(y)$, for some $y$ such that $\phi_{n,\infty}(y)=x^n_j$. Then, the map $\phi_{n,n+1}:\mathfrak{M}_{n+1}\rightarrow \mathfrak{M}_{n}$ between maximal sets induces a map, denoted by: $ \tilde{\phi}_{n,n+1}:\sigma_{B_{n+1}}([D_{n+1},a_{n+1}]) \rightarrow \sigma_{B_n}([D_n,a_n])$ such that:
\begin{equation}
 \tilde{\phi}_{n,n+1}(\ell(a)(x^n_j))=\ell(a)(\phi_{n,n+1}(x^n_j)).
\end{equation} 
Then $\sigma_{B_\omega}([D,a])$ can be identified with the inverse limit given by the inverse system $(\sigma_{B_n}([D_n,a_n]),\mathbb{N}, \tilde{\phi}_{n,n+1})$. Recalling that the manifold $M$ is obtained from the maximal points in $X_\infty$, we deduce that:
\begin{equation}
\sigma_{B_\omega}([D,a])=\left\lbrace \rest{\frac{d a}{d x}}_{x}\right\rbrace_{x\in \mathbb{R}}.
\label{spect}
\end{equation}
Therefore, if we denote by $d_ca$ the de Rham differential of $a$ on $\mathbb{R}$, then we have the following result.
\begin{mprop} (Spectral convergence)
There exists a finite measure $\mu$ and a unitary operator 
\begin{equation}
U:L^2(\mathbb{R})\rightarrow L^2(\mathbb{R},d\mu)
\end{equation}
such that, 
\begin{equation}
U[D,a]U^{-1}\phi=\frac{d a}{d x}\phi, \quad \forall \phi \in L^2(\mathbb{R}),
\end{equation}
Moreover, the norm of $[D,a]$ is given by $\|\left[ D,a\right] \|=\|d_ca\|_\infty$.
\label{spectral}
\end{mprop}
\begin{proof}
This result is an immediate consequence of the spectral theorem on self-adjoint bounded operator (Multiplication operator type) \cite[pp.36-37]{davies_spectral_1995} and the spectrum characterization \eqref{spect}.
\end{proof}
\subsubsection{The d-dimensional lattice}
We start by recalling that the $C^*$-algebra $A(\Lambda)$ is given by the tensor product, 
\begin{equation}
A(\Lambda)=A(L)\otimes \cdots \otimes A(L),
\end{equation}
where $A(L)$ is the algebra associated to the line lattice. Similarly, we consider the sequence of commutative subalgebras 
\begin{equation}
\mathfrak{A}_n(\Lambda)=\mathfrak{A}_n(L)\otimes \cdots \otimes \mathfrak{A}_n(L),
\end{equation}
for every $n\in \mathbb{N}$, that we embed it in the tensor product of matrix algebras
\begin{equation}
B_n=M_{2m_n}(\mathbb{C})\otimes \cdots \otimes M_{2m_n}(\mathbb{C}).
\end{equation}
Again, we adjoin a Dirac operator on each algebra $B_n$ given by:
\begin{equation}
D_n = \sum_{k=1}^{d} 1\otimes \cdots \otimes D^{(k)}_n \otimes \cdots \otimes 1
\label{DiracTensor}
\end{equation}
with the commutator on an element $b$ defined as:
\begin{equation}
[D_n,b]=\sum_{k=1}^{d}b_1\otimes \cdots \otimes [D^{(k)}_n,b_k] \otimes \cdots \otimes b_{d}.
\end{equation}
This gives a spectral triple structure on $B_\omega$ by extending the commutator in the same way as Equation \eqref{commutator}.
\begin{mprop} (Spectral convergence)
There exists a finite measure $\mu$ and a unitary operator 
\begin{equation*}
U:\otimes_{i=1}^d L^2(\mathbb{R})\rightarrow \otimes_{i=1}^d L^2(\mathbb{R},d\mu),
\end{equation*}
such that
\begin{equation*}
U[D,a]U^{-1}\phi=\sum_{k=1}^{d}a_1\phi_1\otimes \cdots \otimes \frac{\partial a_k}{\partial x_k}\phi_k \otimes \cdots \otimes a_{d}\phi_d,
\end{equation*}
for all $\phi=\phi_1\otimes \cdots \phi_k \otimes \cdots \otimes \phi_d$ in $ \otimes_{i=1}^d L^2(\mathbb{R})$.
\end{mprop}
\begin{proof}
Again, we use the fact that the spectrum $[D^{(k)}_n,a_k]$ is given by 
\begin{equation}
\sigma_{B_\omega}([D^{(k)},a_k])=\left\lbrace \rest{\frac{\partial a_k}{\partial x_k}}_{x_k}\right\rbrace_{x_k\in \mathbb{R}},
\end{equation}
for very $k\in \lbrace 1, \dots, d\rbrace$; then using the unitary operator given by
\begin{equation}
U=U_1\otimes \cdots \otimes U_k \otimes \cdots \otimes U_d,
\end{equation}
where for every $k\in \lbrace 1, \dots, d\rbrace$, $U_k$ is the (same) unitary operator given by Proposition \ref{spectral}.
\end{proof}
\subsection{Towards a generalization: the metric question}
\subsubsection{A second example: the torus $\mathbb{T}^d$}
We conclude the present work with the case of the $d$-dimensional torus $\mathbb{T}^d$; we do not specify the metric yet. This example has tow purposes, firstly show how the matrix $D$ depends on the topology of the space $X$ and secondly exhibit the difficulties with respect to the eigenvalues of the commutator. Indeed on the  latter, one has to be able to compute the eigenvalues and moreover, the eigenvalues should somehow reflect the metric (at the limit) the metric $g$ on the manifold $M$.\par 
We start with the case of the circle $S^1$, since the general case of the torus in an arbitrary dimension $d$ is obtained by direct product (similarly to the lattice in $\mathbb{R}^d$). Hence, we are looking for the Dirac operator $D$ associated to a graph obtained from a triangulation of $S^1$. The block matrix $D^{-}$ can be read from Figure \ref{fig:Ncircle} and is given by:
\begin{equation}
D^-=\left( 
\begin{array}{ccccc}
0 & 1 & 0 & \cdots & 0 \\ 
\vdots & \ddots & \ddots & \ddots & \vdots \\ 
\vdots &  & \ddots & \ddots & 0 \\ 
\vdots &  & &\ddots  & 1\\
1 & \cdots & \cdots & \cdots & 0
\end{array} 
\right) .
\end{equation}
where we notice the non-zero coefficients on the down-left corner. Then, we can compute the eigenvalue of the operator $da$. Indeed, for every element $da\in \Omega^1_{D}(\mathfrak{A})$, the spectrum $\sigma (da)$ of the operator $da$ is given by:
\begin{equation}
\sigma(da)=\left\lbrace \pm \frac{1}{h}(\lambda_{j+1} - \lambda_j) : 1\leq j \leq m-1\right\rbrace \cup \left\lbrace  \pm \frac{1}{h}(\lambda_{m} - \lambda_1) \right\rbrace .
\end{equation}
Taking the limit $h\rightarrow 0$, we deduce that there exists a finite measure on $S^1$ and a unitary operator $u$ acting on $L^2(S^1)$ such that
\begin{equation}
\left( u[D,a]u^*\right) \phi=\frac{d a}{d \theta}\phi, \quad \forall \phi \in L^2(S^1).
\end{equation}
Now, we can work out the general case of a $d$-dimensional torus $\mathbb{T}^d$. There exists a unitary operator $u$ acting on $L^2(\mathbb{T}^d)$ such that
\begin{equation*}
\left( u[D,a]u^*\right) \phi=\sum_{k=1}^{d}a_1\phi_1\otimes \cdots \otimes \frac{\partial a_k}{\partial \varphi_k}\phi_k \otimes \cdots \otimes a_{d}\phi_d,
\end{equation*}
for all $\phi=\phi_1\otimes \cdots \phi_k \otimes \cdots \otimes \phi_d$ in $ \otimes_{i=1}^d L^2(S^1)\simeq L^2(\mathbb{T}^d)$.\\

Now, we see that, although there are topological changes in the Dirac operator, the results obtained for the $d$-dimensional torus are very similar, as far as the previous results are concerned, to those obtained for the lattice in $\mathbb{R}^d$. However, seen as Riemannian manifolds, one may expect the approximation of, $\mathbb{R}^d$ with its standard metric on one hand, and of the torus $(\mathbb{T}^d,g)$ with a metric $g$ on the other hand, to reflect somehow the intrinsic geometrical differences. This is one of the main point of the following discussion.
\subsubsection{Discussion}
We notice that the previous results obtained do not depend on the metric $g$ on the torus. Indeed, one could either look at the flat torus with the metric inherited from a quotient of $\mathbb{R}^{d+1}$ or the metric induced from the ambient space it is embedded in. However, as it is defined, there is no use of the metric in $D$; hence, one could, for instance on $S^1$, redefine the matrix $D^-$ as the ansatz
\begin{equation}
D^-=\rho(\theta)\left( 
\begin{array}{ccccc}
0 & 1 & 0 & \cdots & 0 \\ 
\vdots & \ddots & \ddots & \ddots & \vdots \\ 
\vdots &  & \ddots & \ddots & 0 \\ 
\vdots &  & &\ddots  & 1\\
1 & \cdots & \cdots & \cdots & 0
\end{array} 
\right) .
\label{metricDirac1}
\end{equation}
where $\rho$ is a function depending on the metric. Moreover, as a second example, we consider the $2$-dimensional torus $\mathbb{T}^2$ with the parametrization
\begin{equation}
\Psi(\theta,\phi)=\left\langle  \left( R+r\cos\phi\right) \cos \theta,\left( R+r\cos\phi\right)\sin \theta, r\sin \phi\right)
\end{equation}
then inverse of the metric is given (in matrix form) by
\begin{equation}
g^{ij}=(g_{ij})^{-1}=\left(  
\begin{array}{cc}
\frac{1}{R+r\cos \phi} & 0 \\ 
0 & \frac{1}{r^2}
\end{array} 
\right).
\end{equation}
Therefore, we could define the Dirac operator in that case, following Equation \eqref{DiracTensor}, by:
\begin{equation}
D=g^{11}\left( D_1\otimes 1\right)  + g^{22}\left( 1\otimes D_2\right) .
\label{metricDirac2}
\end{equation}
Nevertheless, this can only be a local description i.e. only valid on a chart $(U,\varphi)$ (if the manifold in question is non-trivial). Indeed, one could adopt an \textit{extrinsic} point of view of the Dirac operator: given a local chart $(U,\varphi)$, with a local metric $g|_{U}$, one could consider a lattice-like approximation of $U$ and define a Dirac operator $D|_{U}$ following the same construction than Equations \eqref{metricDirac1} and \eqref{metricDirac2}. However, given two coordinate charts $(U,\varphi_U)$ and $(V,\varphi_V)$ with associated lattices, say $\lambda_U$ and $\lambda_V$ and Dirac operators $D_U$ and $D_V$, it is not clear how to glue them to obtain a lattice $\Lambda_{U\cup V}$ with a Dirac operator $D_{U\cup V}$ that restricts to $D_{U}$ on $\Lambda_U$, respectively $D_{V}$ on $\Lambda_V$.\\

Hence, we can adopt an \textit{intrinsic} construction instead, where the matrix $D$ is defined globally starting from a graph. However, this time we see that the coefficient, say $\omega_{ij}$, of the Dirac operator must depend on the intrinsic geometry of the manifold i.e. the coefficients $\omega_{ij}$ are not simply $0$ or $1$ but are computed from a priori knowledge of the metric $g$.\\

Tied to the previous interrogation on the metric, and the non-triviality of the manifold, is the eigenvalues of a \textit{compatible Dirac operator}. Indeed, in the examples treated above, the correspondence between $D$ and the graph associated to the triangulation gives a commutator for wish the eigenvalues are easy to compute. However,  perform the same task in all generality becomes an intractable problem. Moreover, there is no evidence that given a sequence of commutator, this sequence will converge to the exterior derivative.\\

Therefore, and for all the reasons mentioned above, the fundamental question of the choice of coefficients $\omega_{ij}$ and their relations with the intrinsic properties of the manifold must be tackled in order to see the desired geometry emerge from a sequence of discrete approximations.

\paragraph*{Where are we at now and where are we heading ?} So far, we have defined a spectral triple $(A,H,D)$ on a given triangulation $X$. It is very crucial to notice at this point that one can almost forget about the space $X$ and rely solely on the spectral triple. Indeed, we have shown in Section \ref{Calg} that the algebra $A$ plays the role of functions on $X$ and is enough to recover smooth functions at the limit. Moreover, we have built a correspondence between the given triangulation $X$ and a Dirac operator $D$: the non-zero coefficients of $D$ are determined by the connectivity between vertices of the graph. Hence, $D$ encodes to some extend the topology of $X$. The bracket $[D,a]$ can be then represented as a bounded operator acting on the Hilbert space $H$.\par
However, the last discussion has shown that this is  not enough to represent the metric of the manifold. Thus, we ask now the question on how to set the coefficients $\omega_{ij}$ of $D$ so that at the limit (in the sense of \eqref{commutator}) the sequence converges.
\section{Conclusion}
In this article, we have laid the foundation of a representation theoretic-description of discrete differential calculus using the tools of noncommutative geometry. Starting from a manifold $M$, we construct an inverse system of triangulation, $(K_n)$ which become sufficiently fine for large $n$. We associate to each space $K_n$ a $C^*$-algebra $A_n$ such that the triangulation $K_n$ is identified with its spectrum $Spec(A_n)$. The $C^*$-algebras give a piecewise-linear structure to the triangulations. We then form an inductive system $(A_n)$ with limit a $C^*$-algebra $A_\infty$ with centre isomorphic to the space of continuous function $C(M)$. In this sense, any element $g\in C(M)$ can be uniformly approximated arbitrarily closely by elements $a_n$ in the central subalgebras $\mathfrak{A}_n$. In addition, the sequence of representations $(H_n)$ is also considered as a direct system with limit $H_\infty$ containing the space of square integrable functions $L^2(M)$. Finally, we define spectral triples $(\mathfrak{A},\mathfrak{h},D_n)$ where $D_n$ is the so-called Dirac operator. We show that under certain conditions, the sequence $(D_n)$ converges to the classical Dirac operator in $\mathbb{R}^d$.\par
Our follow-up work focuses on extending the construction presented in this paper beyond the case of the $d$-lattice. We are mainly interested in the Laplacian and Dirac operator in the Riemannian manifold setting.\\
This construction may also provide a unifying framework for geometric discretization such as discrete exterior calculus (DEC) and finite elements exterior calculus (FEEC). The above connection between noncommutative geometry and classical discretization is rather subtle and may lead to powerful and novel numerical approximation.
\bibliographystyle{amsplain}
\bibliography{Article1Bibli.bib}
\end{document}